\documentclass[11pt]{article}

\usepackage{amsmath}
\usepackage{amssymb}

\usepackage{mathtext}
\usepackage[T2A]{fontenc}
\usepackage{latexsym}

\newtheorem{theorem}{Theorem}[section]
\newtheorem{lemma}{Lemma}[section]
\newtheorem{corollary}{Corollary}[section]
\newtheorem{proposition}{Proposition}[section]

\begin{document}

\baselineskip 16pt

\title{\bf  On $\tau$-closed $n$-multiply $\sigma$-local
\\
  formations of finite groups}

\author{Inna N. Safonova \\
{\small Department of Applied Mathematics and Computer Science,}\\
{\small Belarusian State University, Minsk 220030, Belarus}\\ 
{\small E-mail: safonova@bsu.by}}

\date{}
\maketitle

\begin{abstract}  All groups under consideration are finite. Let $\sigma =\{\sigma_{i} \mid i\in I \}$ be some partition of the set of all primes 
$\mathbb{P}$, $G$ be a group, and $\mathfrak F$ be a class of groups. Then 
$\sigma (G)=\{\sigma_i\mid \sigma_i\cap \pi (G)\ne \emptyset\} $ and $\sigma
 (\mathfrak{F})=\cup_{G\in \mathfrak{F}}\sigma (G).$ 
 
A function $f$ of the form
$f:\sigma \to\{\text{formations of groups}\}$ is called
a formation $\sigma$-function. For any formation $\sigma$-function $f$ the class $LF_{\sigma}(f)$ is defined as follows:
$$
LF_{\sigma}(f)=(G \ \text{ is\ a\ group\ } \mid G=1 \ \text{or  }\ G\ne 1\ \text{and \ }\ 
G/O_{\sigma_i', \sigma_i}(G) \in f(\sigma_{i}) \ \text{\ for\ all\ } \sigma_i \in \sigma(G)).
$$ 
If for some formation $\sigma$-function $f$ we have $\mathfrak{F}=LF_{\sigma}(f),$ then the class $\mathfrak{F}$ is called $\sigma $-local and $f$ is called a $\sigma$-local definition of $ \mathfrak{F}.$ 
Every formation is called 0-multiply $\sigma $-local. For $n > 0,$ a formation  $\mathfrak{F}$ is called $n$-multiply $\sigma $-local provided either $\mathfrak{F}=(1)$ is the class of all identity groups 
 or $\mathfrak{F}=LF_{\sigma}(f),$ where $f(\sigma_i)$ is $(n-1)$-multiply $\sigma$-local for all $\sigma_i\in \sigma (\mathfrak{F}).$ A formation is called totally $\sigma$-local if it is $n$-multiply $\sigma$-local for all nonnegative integer $n.$
 
Let $\tau(G)$ be a set of subgroups of $G$ such that 
$G\in \tau(G)$.  Then $\tau$ is called a {subgroup  functor} if for every  epimorphism $\varphi$ : $A \to~B$ and any groups  $H \in \tau (A)$ and $T\in \tau (B)$  we have $H^{\varphi}\in\tau(B)$ and $T^{{\varphi}^{-1}}\in\tau(A)$.
A class  of groups $\mathfrak{F}$ is called 
{$\tau$-closed} if $\tau(G)\subseteq\mathfrak{F}$ for all  $G\in\mathfrak F$.

In this paper,  we describe some properties of $\tau$-closed $n$-multiply $\sigma$-local formations of finite groups, as well as the main properties of the lattice of such formations. In particular, we prove that the set $l^{\tau}_{\sigma_n}$ of all $\tau$-closed $n$-multiply $\sigma$-local formations forms a complete modular algebraic lattice of formations. In addition, we proof that the lattice 
$l^{\tau}_{\sigma_n}$ is $\sigma$-inductive and $\mathfrak G$-separable.
\end{abstract}


\let\thefootnote\thefootnoteorig

\section*{Introduction}

All groups under consideration are finite.  We adhere to the terminology and notation adopted in \cite{sh}--\cite{bbe}. Basic concepts and some properties of $\sigma$-local formations can be found in \cite{skprobl}--\cite{chsafsk}. 

A.N.~Skiba presented \cite{skprobl} the concept of $\sigma$-locality of formations as a tool for studying $\sigma$-properties of groups \cite{skgsu2014, 1}. 
Some interesting applications of $\sigma$-local and $n$-multiply $\sigma$-local formations in the study of finite non-simple groups are received in \cite{skprobl}--\cite{chsafsk}, \cite{alg12}--\cite{skbja}.

In \cite{skprobl}, using $\sigma$-local formations, A.N.~Skiba studied ${\cal S}_{t}^{\sigma}$-closed and ${\cal M}_{t}^{\sigma}$-closed classes of finite groups.  
Some general properties of $\sigma$-local formations as well as their applications for studying $\Sigma_{t}^{\sigma}$-closed classes of meta-$\sigma$-nilpotent groups \cite{zskuk} and ${\Gamma}_{t}^{\sigma}$-closed classes of finite groups \cite{zskac}, were obtained Zhang Chi and A.N.~Skiba. In \cite{chsafskg, chsafsk} Zhang Chi, V.G.~Safonov and A.N.~Skiba  was consider one application of $n$-multiply $\sigma$-local formations in the theory of finite factorisable groups. 
Applications of the theory of $\sigma$-local formations were obtained by A.N.~Skiba \cite{alg12, skbsu} for a characterization of $\sigma$-soluble $P\sigma T$-groups, and also for constructing new sublattices of the lattice of all subgroups of the group generated by formation Fitting sets \cite{skbja}.

Zhang Chi, V.G.~Safonov and A.N.~Skiba \cite{chsafskg, chsafsk} was proved that the lattice $l^{\sigma}_n$ of all $n$-multiply $\sigma$-local formations of finite groups is algebraic and modular.
A.A.~Tsarev \cite{ts} showed that every law of the lattice of all formations is fulfilled in the lattice $l^{\sigma}_n$ 
and that for any nonnegative integer $n$ the lattice $l^{\sigma}_n$ is modular but is not distributive. 
I.N.~Safonova, V.G.~Safonov \cite{safivbsu} was proved that the set $l^\sigma_\infty$ of all totally $\sigma$-local formations of finite groups is a complete algebraic and distributive lattice, and studed also some general properties of totally $\sigma$-local formations of finite groups. I.N.~Safonova \cite{safingsu} studied minimal $\sigma$-local non-$\mathfrak H$-formations, where $\mathfrak H$ is some $\sigma$-local formation. In particular, a description of minimal $\sigma$-local non-$\sigma$-soluble formations of finite groups is obtained.

In this paper, developing the concept of functorial closure of formations proposed by A.N.~Skiba \cite{sk1}, we study the properties of $\tau$-closed $n$-multiply $\sigma$-local formations, as well as the properties of lattices of formations of this type. In particular, we prove that the set $l^{\tau}_{\sigma_n}$ of all $\tau$-closed $n$-multiply $\sigma$-local formations forms a complete modular algebraic lattice of formations. In addition, we proof that the lattice 
$l^{\tau}_{\sigma_n}$ is $\sigma$-inductive and $\mathfrak G$-separable.

\section{Definitions, notations and preliminary statements}

Following \cite{sh}, $\sigma$  is some partition of  
$\mathbb{P}$, that is,  $\sigma =\{\sigma_{i} \mid i\in I \}$, where   $\mathbb{P}=\bigcup_{i\in I} \sigma_{i}$
 and $\sigma_{i}\cap \sigma_{j}= \emptyset  $ for
 all $i\ne j$.  Let $n$ be an integer, $G$ be a group, and $\mathfrak F$ be a class of groups. Then the symbol $\sigma (n)$ denotes the set $\{ \sigma_i\mid \sigma_i \cap \pi (n)\ne \emptyset \},$ $\sigma (G)=\sigma (|G|)$, and $\sigma ({\rm {\mathfrak F}})=\cup _{G\in {\rm {\mathfrak F}}} \sigma (G).$

The group $G$ is called \cite{skgsu2014, 1}: {\sl $\sigma$-primary} if $G$ is a $\sigma_i $-group for some $i;$ {\sl $\sigma$-nilpotent} if every chief factor $H/K$ from $G$ is $\sigma$-central in $G,$ that is, the semidirect product $(H/K){\rm \rtimes }(G/\mbox{C}_{G} (H/K))$ is $\sigma$-primary; {\sl $\sigma$-soluble} if $G=1$ or $G\ne 1$ and each chief factor from $G$ is $\sigma$-primary. The symbol ${\rm {\mathfrak S}}_{\sigma } $ denotes the class of all $\sigma$-soluble groups and ${\rm {\mathfrak N}}_{\sigma } $ denotes the class of all $\sigma$-nilpotent groups. 

Recall that a class of groups ${\rm {\mathfrak F}}$ is called {\sl a formation} if: 1) $G/N\in {\rm {\mathfrak F}}$ when $G\in {\rm {\mathfrak F}},$ and 2)  $G/N\cap R\in {\rm {\mathfrak F}}$ when $G/N\in {\rm {\mathfrak F}}$ and $G/R\in {\rm {\mathfrak F}}.$ If $\mathfrak X$ is a class of groups and $G\in\mathfrak X$, then $G$ is called {\sl a $\mathfrak X$-group.}

Any function $f$ of the form $f:\sigma \to \{ \text{\rm formation of groups}\} $ is called {\sl a formation $\sigma$-function} \cite{skprobl}. For any formation $\sigma$-function $f$ the class $LF_{\sigma } (f)$ defined as follows: 
$$
LF_{\sigma } (f)=(G \; |\; G=1\; \text{\rm or}\; G\ne 1\; \text{\rm and}\; G/O_{\sigma_i ',\sigma_i } (G)\in f(\sigma_i )\; \text{\rm for all\ }\sigma_i \in \sigma (G)).
$$ 

If for some formation $\sigma$-function $f$ the equality ${\rm {\mathfrak F}}=LF_{\sigma } (f),$ then the class ${\rm {\mathfrak F}}$ is called {\sl $\sigma$-local}, and $f$ called {\sl $\sigma$-local definition} of ${\rm {\mathfrak F}}.$ We write $F_{\{ \sigma_i \} } (G)$ instead of $O_{\sigma_i ',\sigma_i } (G)=G_{\mathfrak G_{\sigma_i' }\mathfrak G_{\sigma_i } } .$

Every formation is called {\sl 0-multiply $\sigma $-local} \cite{skprobl}. For $n > 0,$ a formation  $\mathfrak{F}$ is called {\sl $n$-multiply $\sigma $-local} provided either $\mathfrak{F}=(1)$ is the class of all identity groups or $\mathfrak{F}=LF_{\sigma}(f),$ where $f(\sigma_i)$ is $(n-1)$-multiply $\sigma$-local for all $\sigma_i\in \sigma (\mathfrak{F}).$ A formation is called {\sl totally $\sigma$-local} if it is $n$-multiply $\sigma$-local for all nonnegative integer $n.$

Let $\tau(G)$ be a set of subgroups of $G$ such that 
$G\in \tau(G)$.  Then we say, following  \cite{sk1}, 
that $\tau$ is called a {\sl subgroup  functor} if for every  epimorphism
 $\varphi$ : $A \to~B$ and any groups  $H \in \tau (A)$ and $T\in \tau 
(B)$  we have $H^{\varphi}\in\tau(B)$ and $T^{{\varphi}^{-1}}\in\tau(A)$.
A class  of groups $\mathfrak{F}$ is called 
{$\tau$-closed} if $\tau(G)\subseteq\mathfrak{F}$ for all  $G\in\mathfrak F$.

The subgroup functor $\tau$ is called \cite{sk1}: {\sl closed} if for any two groups $G$ and $H\in\tau(G)$ we have $\tau(H)\subseteq\tau(G)$; {\sl trivial}, if for any group $G$ we have $\tau(G) = \{G\}$;  {\sl single}, if for any group $G$ we have $\tau(G)$  is the collection of all subgroups of $G$.

A class of groups $\mathfrak F$ is called {\sl $\tau$-closed}, if $\tau(G)\subseteq\mathfrak F$ for any group $G\in\mathfrak F$. A formation is called {\sl $\tau$-closed}, if it is a $\tau$-closed class of groups.

The collection of all $\tau$-closed $n$-multiply $\sigma$-local formations denote by $l^{\tau}_{\sigma_n}$. 
Formations from $l^{\tau}_{\sigma_n}$ we call {\sl $l^{\tau}_{\sigma_n}$-formations.}
For any set of groups ${\rm {\mathfrak X}}$ the symbol $l^{\tau}_{\sigma_n} {\rm form\,}({\rm {\mathfrak X}})$ denotes {\sl a $\tau$-closed $n$-multiply $\sigma$-local formation generated by ${\rm {\mathfrak X}},$} that is, $l^{\tau}_{\sigma_n} {\rm form\,}({\rm {\mathfrak X}})$ is the intersection of all $\tau$-closed $n$-multiply $\sigma$-local formations containing ${\rm {\mathfrak X}}.$ If ${\rm {\mathfrak X}}=\{ G\} $ for some group $G,$ then ${\rm {\mathfrak F}}=l^{\tau}_{\sigma_n} {\rm form\,}(G)$ is called {\sl a one-generated} $\tau$-closed $n$-multiply $\sigma$-local formation. Let  $\{ \mathfrak{F}_j \mid j\in J \}$ be some collection of $\tau$-closed $n$-multiply $\sigma$-local formations. Then we put
$$
\vee ^{\tau}_{\sigma_n}(\mathfrak{F}_j \mid j\in J )=l^{\tau}_{\sigma_n} {\rm form\,}( \cup_ {j\in J}\mathfrak{F}_j).
$$
In particular, for any two $l^{\tau}_{\sigma_n} $-formations $\mathfrak M$ and ${\rm {\mathfrak H}}$ we set $\mathfrak M\vee ^{\tau}_{\sigma_n}\mathfrak H=l^{\tau}_{\sigma_n} {\rm form\,}( \mathfrak M\cup \mathfrak H).$

If $f$ is a formation $\sigma$-function, then the symbol ${\rm Supp}(f)$ denotes the support of $f,$ that is, the set of all $\sigma_i $ such that $f(\sigma_i )\ne \emptyset .$ 
A formation $\sigma$-function $f$ is called: {\sl $l^{\tau}_{\sigma_n}$-valued}, if $f(\sigma_i)$ is $\tau$-closed $n$-multiply $\sigma$-local formation for each $\sigma_i\in {\rm Supp}(f)$; {\sl integrated} if $f(\sigma_i )\subseteq LF_{\sigma } (f)$ for all $i$; {\sl full} if $f(\sigma_i)=\mathfrak G_{\sigma_i}f(\sigma_i)$ for all $i$. 

If  $F$ is a full integrated formation $\sigma$-function and
 $\mathfrak{F}=LF_{\sigma}(F)$, then $F$ is called the {\sl canonical $\sigma$-local definition} of $\mathfrak{F}$. 

Following \cite{sh}, on the set of all $l^{\tau}_{\sigma_n}$-valued formation $\sigma$-functions we define the partial order $''\le \ ''$ as follows. Let  $f_1$ and $f_2$ be $l^{\tau}_{\sigma_n}$-valued formation $\sigma$-functions, then $f_1\le f_2$ if and only if $f_1(\sigma_i)\subseteq f_2(\sigma_i)$ for all $\sigma_i\in \sigma.$

Let $\{f_j\mid j\in J\}$ be some collection of $l^{\tau}_{\sigma_n}$-valued formation $\sigma$-functions. 
We define the $l^{\tau}_{\sigma_n}$-valued formation $\sigma$-function $f=\vee^{\tau}_{\sigma_n}(f_j\mid j\in J)$ as follows:
$$
f(\sigma_i)=\vee^{\tau}_{\sigma_n}(f_j\mid j\in J)(\sigma_i)=\vee^{\tau}_{\sigma_n}(f_j(\sigma_i)\mid j\in J)=l^{\tau}_{\sigma_n} {\rm form\,}(\cup_{j\in J} f_j(\sigma_i)), 
$$ 
in particular, 
$$
f(\sigma_i)=(f_1\vee^{\tau}_{\sigma_n} f_2)(\sigma_i)=f_1(\sigma_i)\vee^{\tau}_{\sigma_n} f_2(\sigma_i)=l^{\tau}_{\sigma_n} {\rm form\,}(f_1(\sigma_i)\cup f_2(\sigma_i)),
$$
if at least one of the formations $f_j(\sigma_i)\ne\emptyset$ and $f(\sigma_i)=\emptyset$ if $f_j(\sigma_i)=\emptyset$ for all $j\in J$. 

We also use $\cap_{j\in J}f_j$ to denote a formation $\sigma$-function $h$ such that $h(\sigma_i)=\cap_{j\in J}f_j(\sigma_i)$, in 
particular, $h(\sigma_i)=(f_1\cap f_2)(\sigma_i )=f_1(\sigma_i)\cap f_2(\sigma_i)$, for all $i.$
 
Let $\{f_j\mid j\in J\}$ be a collection of all $l^{\tau}_{\sigma_n}$-valued $\sigma$-local definitions of $\mathfrak F.$ Then we say that $f=\cap_{j\in J}f_j$ is {\sl the smallest $l^{\tau}_{\sigma_n}$-valued $\sigma$-local definition} of $\mathfrak F.$

Recall \cite{chsafsk} that for an arbitrary collection of groups $\mathfrak{X}$ and any $\sigma_i\in\sigma$ by the symbol $\mathfrak{X}(\sigma_i)$ denote class of groups, defined as follows:

$$
\mathfrak{X}(\sigma_i) =\begin{cases}
(G/F_{\{\sigma_i\}}(G)\mid G \in \mathfrak{X}), & \mbox{ if \  } \sigma_i\in\sigma (\mathfrak{X}), \cr
\emptyset, & \mbox{ if \  } \sigma_i\notin\sigma (\mathfrak{X}).
\end{cases}
$$

We will need a number of known results, which we formulate in the form of the following lemmas

\begin{lemma} {\rm \cite{safivbsu}}\label{l1}
If $\mathfrak F=LF_\sigma(f)$ and $G/\mbox{O}_{\sigma_i}(G)\in f(\sigma_i)\cap \mathfrak F$ for some $\sigma_i~\in~\sigma(G),$ then $G\in\mathfrak F.$
\end{lemma}

\begin{lemma} {\rm \cite{zskuk}}\label{l2} 
{\rm (1)} For every formation $\sigma$-function $f$ the class $LF_{\sigma}(f)$ is a nonempty saturated formation.

{\rm (2)} Every $\sigma$-local formation $\mathfrak F$ has an unique $\sigma$-local definition $F$ such that, for any $\sigma$-local definition $f$ of the formation $\mathfrak F$ and any $\sigma_i\in\sigma$ the following relation is true:
$$
F(\sigma_i)=\mathfrak G_{\sigma_i}(f(\sigma_i)\cap \mathfrak F) =\mathfrak G_{\sigma_i}F(\sigma_i).
$$ 
\end{lemma}

\begin{lemma} {\rm \cite{zskuk}}  \label{l3}
If $\mathfrak{F}=LF_{\sigma}(f)$, then $\mathfrak{F}=LF_{\sigma}(t)$, where $t(\sigma_i)=f(\sigma_i)\cap \mathfrak{F}$ for all $\sigma_i\in\sigma$.
\end{lemma}

\begin{lemma} {\rm \cite{chsafsk}} \label{l4} 
If $\mathfrak{F}=\cap _{j\in J}\mathfrak{F}_{j}$ and $\mathfrak{F}_{j}=LF_{\sigma}(f_{j})$ for all $j\in J$, then $\mathfrak{F}=LF_{\sigma}(f)$, where $f(\sigma_i)=\cap _{j\in J}f_{j}(\sigma_i)$ for all $\sigma_i\in \sigma (\mathfrak{F})=\cap _{j\in J}\sigma(\mathfrak{F}_{j})$ and $f(\sigma_i)=\varnothing$ for all $\sigma_i\in \sigma \setminus\sigma (\mathfrak{F})$. Moreover, if $f_{j}$ is an integrated formation $\sigma$-function for all $j\in J$, then $f$ is also an integrated formation $\sigma$-function.
\end{lemma}  

\begin{lemma} {\rm \cite{chsafsk}}  \label{l5}
If the class of groups ${\rm {\mathfrak F}}_{j} $ is an $n$-multiply $\sigma$-local formation for all $j\in J,$ then the class $\cap _{j\in J} {\rm {\mathfrak F}}_{j} $ is also an $n$-multiply $\sigma$-local formation.
\end{lemma}

\begin{lemma} \label{l6} 
 {\rm \cite{zskuk}}
Let $f$ and $h$ be formation $\sigma$-functions and let $\Pi=\text{\rm Supp}(f)$. Suppose $\mathfrak{F}=LF_{\sigma}(f)=LF_{\sigma}(h)$: 

{\rm (1)}  $\Pi=\sigma (\mathfrak{F})$.  
 
{\rm (2)}  $\mathfrak{F}= (\cap _{\sigma_i\in \Pi}\mathfrak{G}_{\sigma_{i}'} \mathfrak{G}_{\sigma_i}f(\sigma_i)) \cap \mathfrak{G}_{\Pi }.$ Therefore, $\mathfrak{F}$ is a saturated formation.

{\rm (3)} If every group from $\mathfrak{F}$ $\sigma$-soluble, then $\mathfrak{F}= (\cap _{\sigma_i\in \Pi}\mathfrak{S}_{\sigma_{i}'} \mathfrak{S}_{\sigma_i}f(\sigma_i)) \cap  \mathfrak{S}_{\Pi }.$
 
{\rm (4)} If $\sigma_i\in \Pi $, then $\mathfrak{G}_{\sigma_i}(f(\sigma_i)\cap \mathfrak{F})=\mathfrak{G}_{\sigma_i}(h(\sigma_i)\cap \mathfrak{F})\subseteq \mathfrak{F}.$

{\rm (5)} $\mathfrak{F}=LF_{\sigma}(F)$, where $F$ is the single formation $\sigma$-function such that $F(\sigma_i)=\mathfrak{G}_{\sigma_i}F(\sigma_i)\subseteq \mathfrak{F}$ for all $\sigma_i\in \Pi$ and $F(\sigma_i)=\varnothing$ for all $\sigma_i\in \Pi'$. Besides, $F(\sigma_i)=\mathfrak{G}_{\sigma_i}(f(\sigma_i)\cap  \mathfrak{F})$ for all $i$.
\end{lemma} 

\begin{lemma} {\rm \cite{chsafsk}}\label{l7} 
{\rm (i)} Let $ \mathfrak{F}=LF_{\sigma}(f)$, where $f(\sigma_i)=\mathfrak{F}$ for all $\sigma_i\in \sigma (\mathfrak{F})$. Then $ \mathfrak{F}$ is a totally $\sigma$-local formation.

{\rm (ii)} Let $\Pi\subseteq\sigma$ and let $f$ be a formation $\sigma$-function, defined as follows:
$$ f(\sigma_i)=\begin{cases} {\mathfrak{G}}_{\Pi} &\text{if \  $\sigma_i\in \Pi$},\\ 
\varnothing,&\text{if \  $\sigma_i\in \Pi'$}.
\end{cases}
$$ 
Then $LF_{\sigma}(f)= \mathfrak{G}_{\Pi }$ is a totally $\sigma$-local formation. In particular, the formation ${\mathfrak{S}}_{\Pi }$ is totally $\sigma$-local.
\end{lemma} 

\begin{lemma} {\rm \cite{svsib}} \label{l8} 
Let $\mathfrak M$~be a nonempty hereditary formation, $\mathfrak F$~a nonempty $\tau$-closed formation. Then $\mathfrak {MF}$~is a $\tau$-closed formation.
\end{lemma}

\begin{lemma}{\rm \cite{chsafsk}} \label{l9} 
{The set ${\cal S}^{\sigma}_{n}$, of all $n$-multiply $\sigma $-local formations, forms a subsemigroup of the semigroup of all formations $G{\mathfrak{G}}$. Moreover, $|{\cal S}^{\sigma}_{n}|=2^{\aleph _{0}}$ for every $\sigma $ for $|\sigma| > 1$, and ${\mathfrak{G}}_{\sigma_i}$ is the minimal idempotent in ${\cal S}^{\sigma}_{n}$ for all $n > 0$ and $i\in I$.}
\end{lemma}

\begin{lemma} {\rm \cite{chsafsk}} \label{l10} 
If ${\rm {\mathfrak F}}$ is a nonempty formation and $f(\sigma_i )={\rm {\mathfrak F}}$ for all $i,$ then $LF_{\sigma } (f)={\rm {\mathfrak N}}_{\sigma } {\rm {\mathfrak F}}.$
\end{lemma}

\begin{lemma}{\rm \cite[p.~25]{sk1}} \label{l13} 
Let $\mathfrak{X}$ be a $\tau$-closed semiformation and $A\in\mathfrak{F} = \tau\mbox{\rm form\,}\mathfrak{X}$. Then if $A$ is a monolithic group and $A\notin\mathfrak{X}$, then $\mathfrak{F}$ contains a group $H$ with normal subgroups $N, M, N_1,..., N_t,M_1,..., M_t \ (t\ge 2)$, such that the following statements are true:

1) $H/N\simeq A$, $M/N=\mbox{Soc}(H/N)$;

2) $N_1\cap\ldots \cap N_t=1;$

3) $H/N_i$ is a monolithic $\mathfrak{X}$-group with monolith $M_i/N_i$, which is $H$-isomorphic to $M/N$;

4) $M_1\cap \ldots  \cap M_t\subseteq M$.
\end{lemma}

Recall that for any collection of groups $\mathfrak{X}$ by the symbol S$_{\tau}\mathfrak{X}$ \cite[p. 23]{sk1} denote the set of all groups $H$, such that $H\in\tau(G)$ for some group $G\in\mathfrak{X}$ and a class of groups $\mathfrak{F}$ is called {\sl $\tau$-closed}, if S$_\tau(\mathfrak{F}) = \mathfrak{F}$. 

A group $G\in \textnormal{\small{R}}_{0}\mathfrak X $ \cite[p. 264]{dh} if and only if $G$ has normal subgroups $N_1, \ldots, N_t$ $(t\ge 2)$ such that $G/N_i\in \mathfrak X$ and $\cap_{i=1}^t N_i=1$. 
The set of all homomorphic images of all groups from $\mathfrak{X}$ is denoted by Q$\mathfrak{X}$. A class of groups $\mathfrak{F}$ is called {\sl a semiformation} if $\mathfrak{F}$ = Q$\mathfrak{F}$. 
The intersection of all $\tau$-closed semiformations, which contain a given set of groups $\mathfrak{X}$, is called {\sl a $\tau$-closed semiformation generated by $\mathfrak{X}$} \cite[p. 23]{sk1}.

Let $\tau$ be an arbitrary subgroup functor, $\overline {\tau}$ be the intersection of all closed subgroup functors $\tau_i$, for which $\tau \le \tau_i$. The functor $\overline {\tau}$ is called {\sl a closure}  of the functor $\tau$ \cite[p. 20]{sk1}.

\begin{lemma} {\rm \cite[p. 23]{sk1}}\label{l14} 
Let $\mathfrak{F}$ be the $\tau$-closed semiformation generated by $\mathfrak{X}$. Then
$
\mathfrak{F} = {\mbox{\rm QS}}_{\overline\tau}\,\mathfrak{X}.
$
\end{lemma}

\begin{lemma} {\rm \cite[p. 24]{sk1}}\label{l15} 
For any collection of groups $\mathfrak{X}$ we have
$
\tau\mbox{\rm form\,}\mathfrak{X} = \mbox{\rm QR}_{0}\mbox{\rm S}_{\overline\tau}(\mathfrak{X}).
$
\end{lemma}

\begin{lemma} {\rm \cite[p. 24]{sk1}}\label{l16} 
For any collection of $\tau$-closed formations $\{\mathfrak{M}_i\mid i\in I\}$ we have
$$
\tau\mbox{\rm form\,}(\cup_{i\in I}\mathfrak{M}_i) = \mbox{\rm form\,}(\cup_{i\in I}\mathfrak{M}_i).
$$
\end{lemma}

\section{\bf Properties of $l^{\tau}_{\sigma_n}$-formations}

\begin{lemma} \label{l16a}
Let $\mathfrak M=LF_\sigma(m)$, where $m$ is a formation $\sigma$-function such that $m(\sigma_i)$ is a $\tau$-closed formation for all $\sigma_i\in\text{\rm Supp}(m)$. Then $\mathfrak M$ is $\tau$-closed. 
\end{lemma}

Proof.  Let $G\in\mathfrak M$ and $H\in\tau(G)$. Let us show that $H\in\mathfrak M.$ Note that $G/F_{\{\sigma_i\}}(G)\in m(\sigma_i)$ for all $\sigma_i\in\sigma(G)$ since $G\in \mathfrak M$. Hence  $\tau(G/F_{\{\sigma_i\}}(G))\subseteq m(\sigma_i)$ for all $\sigma_i\in\sigma(G)$, since $m(\sigma_i)$ is $\tau$-closed by the condition of the lemma.  

Now we show that $H/F_{\{\sigma_i\}}(H)\in m(\sigma_i)$ for all $\sigma_i\in\sigma(H)$. Let $\sigma_i\in\sigma(H)$. Then, obviously, $\sigma_i\in\sigma(G)\subseteq\text{\rm Supp}(m)$. 
If $H\subseteq F_{\{\sigma_i\}}(G)$, then 
$H=F_{\{\sigma_i\}}(H)$ and $H/F_{\{\sigma_i\}}(H)\simeq 1\in m(\sigma_i)\ne\emptyset$. Let $H\not \subseteq F_{\{\sigma_i\}}(G)$.  Then since $HF_{\{\sigma_i\}}(G)/F_{\{\sigma_i\}}(G)\in\tau(G/F_{\{\sigma_i\}}(G))$, we have 
$$
H/F_{\{\sigma_i\}}(G)\cap H\simeq HF_{\{\sigma_i\}}(G)/F_{\{\sigma_i\}}(G)\in m(\sigma_i).
$$ 
But $F_{\{\sigma_i\}}(G)\cap H\subseteq F_{\{\sigma_i\}}(H)$. Hence  
$$
H/F_{\{\sigma_i\}}(H)\simeq
(H/F_{\{\sigma_i\}}(G)\cap H)/(F_{\{\sigma_i\}}(H)/F_{\{\sigma_i\}}(G)\cap H) \in m(\sigma_i).
$$
Thus, $H/F_{\{\sigma_i\}}(H) \in m(\sigma_i)$ for all $\sigma_i\in\sigma(H)$,  so $H\in \mathfrak M$. Consequently, $\mathfrak M$ is a $\tau$-closed formation. The lemma is proved.

Following \cite[p. 31]{sk1} for any non-negative integer $n$ we define the class of groups
$$
\mathfrak{X}_{\sigma_n}^{\tau}(\sigma_i)=\begin{cases}
l_{\sigma_n}^{\tau}\mbox{form\,}(\mathfrak{X}(\sigma_i)),& \mbox{ if \ } \sigma_i\in\sigma(\mathfrak{X}); \cr \emptyset, & \mbox{ if \ } \sigma_i\in\sigma\setminus\sigma(\mathfrak{X}). \cr
\end{cases}
$$
In particular, if $n=0$, then
$$
\mathfrak{X}_{\sigma_0}^{\tau}(\sigma_i)=\begin{cases} \tau\mbox{form\,}(\mathfrak{X}(\sigma_i)),& \mbox{ if \ } \sigma_i\in\sigma(\mathfrak{X}); \cr \emptyset, & \mbox{ if \ } \sigma_i\in\sigma\setminus\sigma(\mathfrak{X}). \cr
\end{cases}
$$

\begin{theorem}\label{t1} 
Let $\mathfrak{X}$ be some nonempty collection of groups, $\mathfrak{F}=l^{\tau}_{\sigma_n}\text{\rm form\,}\mathfrak{X}=LF_{\sigma}(f)$, where $f$ is the smallest $l^{\tau}_{\sigma_{n-1}}$-valued $\sigma$-local definition of $\mathfrak{F}$. Then the following hold: 

(1) $\sigma(\mathfrak{X}) = \sigma(\mathfrak{F})$;

(2) $f(\sigma_i)=\mathfrak{X}_{\sigma_n}^{\tau}(\sigma_i)=\mathfrak{F}_{\sigma_n}^{\tau}(\sigma_i)$ for all $\sigma_i\in\sigma(\mathfrak{X})$ and $f(\sigma_i)=\emptyset$ for all $\sigma_i\in \sigma\setminus \sigma(\mathfrak{X})$;

(3) if $h$ is an arbitrary $l^{\tau}_{\sigma_{n-1}}$-valued $\sigma$-local definition of $\mathfrak{F}$, then for all $\sigma_i\in \sigma(\mathfrak{X})$ we have 
$$
f(\sigma_i) = l^{\tau}_{\sigma_{n-1}}\text{\rm form\,}(A\mid A\in \mathfrak{F}\cap h(\sigma_i), \mbox{O}_{\sigma_i}(A)=1).
$$

\end{theorem}

Proof. First we show that statement (2) is true. Let $m$ be a formation $\sigma$-function such that $m(\sigma_i) = l^{\tau}_{\sigma_{n-1}}\text{\rm form\,}(\mathfrak{X}(\sigma_i))$ for all $\sigma_i\in\sigma(\mathfrak{X})$ and $m(\sigma_i)=\emptyset$ for all $\sigma_i\in \sigma\setminus \sigma(\mathfrak{X})$.
Let us show that $m=f$. Let $\mathfrak{M}= LF_\sigma(m)$. Then $\mathfrak M$ is $n$-multiply $\sigma$-local by definition. In view of Lemma~\ref{l16a} the formation $\mathfrak M$ is $\tau$-closed. Hence $\mathfrak M\in l^{\tau}_{\sigma_{n}}$. Since  $A/F_{\{\sigma_i\}}(A)\in m(\sigma_i)$ for any group $A\in\mathfrak X$ and for all $\sigma_i\in\sigma(A)$, we have $\mathfrak{X}\subseteq\mathfrak{M}$.
Hence, $\mathfrak{F}\subseteq\mathfrak{M}$. Let $f_1$ be an arbitrary $l^{\tau}_{\sigma_{n-1}}$-valued $\sigma$-local definition of $\mathfrak{F}$.
Then, for any $\sigma_i\in\sigma$ we have $\mathfrak{X}(\sigma_i)\subseteq f_1(\sigma_i)$ since $\mathfrak{X}\subseteq\mathfrak{F}$.
Hence, $m(\sigma_i)\subseteq f_1(\sigma_i)$. Therefore,  $\mathfrak{M}\subseteq\mathfrak{F}$ and
$\mathfrak{M}=\mathfrak{F}$, so $m = f$. Thus, statement (2) is true. In addition, it is obvious that the validity of statement (2) implies the validity of (1).

Let us now prove (3). Let $f_1$ be a formation $\sigma$-function such that
$$
f_1(\sigma_i) = l^{\tau}_{\sigma_{n-1}}\mbox{\rm form\,}(A\mid A\in \mathfrak{F}\cap h(\sigma_i), \mbox{O}_{\sigma_i}(A)=1)
$$
for all $\sigma_i\in\sigma$ and $f_1(\sigma_i) = \emptyset$ for all $\sigma_i\in \sigma\setminus \sigma(\mathfrak{X})$.
Since $f\le h$ and $\mbox{O}_{\sigma_i}(G/F_{\{\sigma_i\}}(G))=1$ for any group $G$ and every $\sigma_i\in\sigma$, we have $f\le f_1$ in view of (2). 
Now let $\sigma_i\in \sigma(\mathfrak{X}),$ $A\in \mathfrak{F}\cap h(\sigma_i)$ and O$_{\sigma_i}(A)=1$.
Let $P$ be a non-identity $\sigma_i$-group and denote by $G$ a regular wreath product $P\wr A=K\rtimes A$, where $K$ is the base group of 
$G$.
Then since $\mbox{O}_{\sigma_i'}(G)=1$, we have $F_{\{\sigma_i\}}(G)=\mbox{O}_{\sigma_i}(G)$ and 
$$\mbox{O}_{\sigma_i}(G)=\mbox{O}_{\sigma_i}(G)\cap KA=K(\mbox{O}_{\sigma_i}(G)\cap A)=K.
$$ 
Therefore, $A \simeq G/K= G/\mbox{O}_{\sigma_i}(G)=G/F_{\{\sigma_i\}}(G)$. 

Since $G/\mbox{O}_{\sigma_i}(G)\simeq A\in \mathfrak{F}\cap h(\sigma_i)$, $G\in \mathfrak{F}$ by Lemma~\ref{l1}.
Hence, $A\simeq G/F_{\{\sigma_i\}}(G)\in f(\sigma_i)$. Therefore,  $f_1(\sigma_i)\le f(\sigma_i)$, so $f=f_1$. The theorem is proved.

\begin{corollary}\label{c1t1} 
Let $f_j$ be the smallest $l^{\tau}_{\sigma_{n-1}}$-valued $\sigma$-local definition
of $\mathfrak{F}_j,$  $j = 1,2$.
Then $\mathfrak{F}_1\subseteq\mathfrak{F}_2$ if and only if $f_1 \le f_2$.
\end{corollary}

In particular, if $\tau$ is a trivial subgroup functor from Theorem~\ref{t1}  is obtained

\begin{corollary}\label{c2t1} 
Let $\mathfrak{X}$ be some nonempty collection of groups, $\mathfrak{F}=l^{\sigma}_n\text{\rm form\,}\mathfrak{X}=LF_{\sigma}(f)$, where  $f$ is the smallest $l^{\sigma}_{n-1}$-valued $\sigma$-local definition of $\mathfrak{F}$. Then the following hold: 

(1) $\sigma(\mathfrak{X}) = \sigma(\mathfrak{F})$;

(2) $f(\sigma_i)=\mathfrak{X}^{\sigma}_{n-1}(\sigma_i)=\mathfrak{F}^{\sigma}_{n-1}(\sigma_i)$ for all $\sigma_i\in\sigma(\mathfrak{X})$ and $f(\sigma_i)=\emptyset$ for all $\sigma_i\in \sigma\setminus \sigma(\mathfrak{X})$;

(3) if $h$ is an arbitrary $l^{\sigma}_{n-1}$-valued $\sigma$-local definition of $\mathfrak{F}$, then for all $\sigma_i\in \sigma(\mathfrak{X})$ we have 
$$
f(\sigma_i) = l^{\sigma}_{n-1}\text{\rm form\,}(A\mid A\in \mathfrak{F}\cap h(\sigma_i), \mbox{O}_{\sigma_i}(A)=1).
$$
\end{corollary}

\begin{corollary}\label{c3t1}  
Let $f_j$ be the smallest $l^{\sigma}_{n-1}$-valued $\sigma$-local definition of $\mathfrak{F}_j,$  $j = 1,2$. Then $\mathfrak{F}_1\subseteq\mathfrak{F}_2$ if and only if  $f_1 \le f_2$.
\end{corollary}

\begin{corollary}\label{c4t1} 
Let $\mathfrak{X}$ be a nonempty class of groups, $\mathfrak{F}= l_{\sigma}\text{\rm form\,}\mathfrak{X}=LF_{\sigma}(f)$, where $f$ is the smallest $\sigma$-local definition of $\mathfrak{F}$. Then the following hold: 

1) $\sigma(\mathfrak{X}) = \sigma(\mathfrak{F})$;

2)  $f(\sigma_i)=\text{\rm form\,}(\mathfrak{X}(\sigma_i))$ for all $\sigma_i\in\sigma$;

3) if $h$ is an arbitrary $\sigma$-local definition of $\mathfrak{F}$, then   $f(\sigma_i)=\text{\rm form\,}(A\mid A\in \mathfrak{F}\cap h(\sigma_i), \mbox{O}_{\sigma_i}(A)=1)$ for all $\sigma_i\in \sigma(\mathfrak{F})$.
\end{corollary}

\begin{corollary} \label{c5t1} 
{\rm \cite{safingsu}} 
Let $f_j$ be the smallest $\sigma$-local definition of $\mathfrak{F}_j,$  $j = 1,2$. Then $\mathfrak{F}_1\subseteq\mathfrak{F}_2$ if and only if $f_1 \le f_2$.
\end{corollary}

Recall what if $\mathfrak{X}$ is some a nonempty collection of groups, $p\in\mathbb{P}$, then the classes $\mathfrak{X}(p)$ and $\mathfrak{X}^\tau_n(p)$ are defined as follows \cite[p. 31]{sk1}: 
$$
\mathfrak{X}(p) =\begin{cases}
\mbox{form\,}(G/\mbox{F}_p(G)\mid G \in \mathfrak{X}), & \mbox{ if \ } p\in\pi (\mathfrak{X}), \cr
\emptyset, & \mbox{ if \ } p\notin\pi (\mathfrak{X}),
\end{cases}
$$
$$
\mathfrak{X}^\tau_n(p) =\begin{cases}
l^{\tau}_n\mbox{form\,}(\mathfrak{X}(p)), & \mbox{ if \ } p\in\pi (\mathfrak{X}), \cr
\emptyset, & \mbox{ if \ } p\notin\pi (\mathfrak{X}).
\end{cases}
$$

In the classical case, when $\sigma=\sigma^1=\{ \{2\}, \{3\}, \ldots \}$ from Theorem~\ref{t1} we obtain the following well-known results

\begin{corollary} {\rm\cite[p. 32]{sk1}} \label{c6t1} 
Let $\mathfrak{F}=l^{\tau}_n\text{\rm form\,}\mathfrak{X}$, where $\mathfrak{X}$ is a nonempty group class, $n\ge 1.$ If $f$ is the minimal $l^{\tau}_{n-1}$-valued screen of $\mathfrak{F}$, then the following hold: 

1) $\pi(\mathfrak{X}) = \pi(\mathfrak{F})$;

2) $f(p)=\mathfrak{X}^\tau_n(p)=\mathfrak{F}^\tau_n(p)$ for all $p\in\pi(\mathfrak{X})$;

3) if $h$ is an arbitrary $l^{\tau}_{n-1}$-valued screen of $\mathfrak{F}$, then for all $p\in \pi(\mathfrak{X})$ we have 
$$
f(p) = l^{\tau}_{n-1}\text{\rm form\,}(A\mid A\in \mathfrak{F}\cap h(p), O_{p}(A)=1).
$$
\end{corollary}

\begin{corollary} \label{c7t1} 
{\rm \cite[p. 30]{sk1}} 
Let $f_j$ be the minimal $l^{\tau}_{n-1}$-valued screen of $\mathfrak{F}_j,$  $j = 1,2$.
Then $\mathfrak{F}_1\subseteq\mathfrak{F}_2$ if and only if $f_1 \le f_2$.
\end{corollary}

\begin{corollary} {\rm\cite[p. 15]{sk1}}\label{c8t1} 
Let $\mathfrak{X}$ be nonempty class of groups, $\mathfrak{F}=l\text{\rm form\,}\mathfrak{X}$ and $f$ be the minimal local screen of $\mathfrak{F}$. Then the following hold: 

1) $\pi(\mathfrak{X}) = \pi(\mathfrak{F})$;

2)  $f(p)=\mathfrak{X}(p)$ for all $p\in\mathbb{P}$;

3) if $h$ is an arbitrary local screen of $\mathfrak{F}$, then  $ f(p) = \text{\rm form\,}(A\mid A\in \mathfrak{F}\cap h(p), O_{p}(A)=1)$ for all $p\in \pi(\mathfrak{F})$.
\end{corollary}

\begin{corollary} \label{c9t1} 
{\rm \cite[p. 15]{sk1}} 
Let $f_j$ be the minimal local screen of $\mathfrak{F}_j,$  $j = 1,2$. Then $\mathfrak{F}_1\subseteq\mathfrak{F}_2$ if and only if $f_1 \le f_2$.
\end{corollary}

Following \cite[p. 27]{sk1}, we introduce the concept of the $\sigma$-locality index of a formation. Let $n$ be a non-negative integer. We say that the formation $\mathfrak{F}$ has a {\sl $\sigma$-locality index} $\text{\rm Ind}_\sigma(\mathfrak{F})$ equal to $n$, if $\mathfrak{F}$ is a $n$-multiply $\sigma$-local formation, but is not $(n+1)$-multiply $\sigma$-local. 
If the formation $\mathfrak{F}$ is $n$-multiply $\sigma$-local for all natural numbers $n$, then we say that $\mathfrak{F}$ has {\sl infinite index $\sigma$-locality} and write $\text{\rm Ind}_\sigma(\mathfrak{F})=\infty$. If the formation $\mathfrak{F}$ is not $\sigma$-local, then $\text{\rm Ind}_\sigma(\mathfrak{F})=0$.

It is clear that for a totally $\sigma$-local formation $\mathfrak F$ we have $\text{\rm Ind}_\sigma(\mathfrak{F})=\infty$. In particular, as you know \cite{chsafsk}, the formation (1) of all identity groups, the formation of $\mathfrak N_{\sigma}$ of all $\sigma$-nilpotent groups, and the formation $\mathfrak S_{\sigma}$ of all $\sigma$-soluble groups are totally $\sigma$-local. Therefore $\text{\rm Ind}_\sigma((1))=\infty$, $\text{\rm Ind}_\sigma(\mathfrak N_{\sigma})=\infty$ and $\text{\rm Ind}_\sigma(\mathfrak S_{\sigma})=\infty$. 

In particular, in the classical case, when $\sigma =\sigma ^{1} =\{ \{ 2\} ,\{ 3\} ,\ldots \} $ a $\sigma^1$-locality index of a formation, is obviously, coincides with the formation locality indix, i.e. $\text{\rm Ind}_{\sigma^1}(\mathfrak{F})=\text{\rm Ind}_l(\mathfrak{F})$.

\begin{proposition}\label{p1}  
Let $F$ be the canonical $\sigma$-local definition of $\mathfrak{F}$.
Then $\text{\rm Ind}_\sigma(\mathfrak{F})=n$ if and only if there is $\sigma_i\in\sigma$ such that $\text{\rm Ind}_\sigma(F(\sigma_i))=n-1$, and $\text{\rm Ind}_\sigma(F(\sigma_j))\ge n-1$ for all $\sigma_j\in\sigma(\mathfrak{F})\setminus\{\sigma_i\}$.
\end{proposition}

Proof. Necessity. Let $\text{\rm Ind}_\sigma(\mathfrak{F})=n$. Then  $\mathfrak{F}$ has a $\sigma$-local definition $f$, such that $\text{\rm Ind}_\sigma(f(\sigma_i))=n-1$ for some $\sigma_i\in\sigma(\mathfrak{F})$ and $\text{\rm Ind}_\sigma(f(\sigma_j))\ge n-1$ for all $\sigma_j\in\sigma(\mathfrak{F})\setminus\{\sigma_i\}$.
In addition, for any other $\sigma$-local definition $k$ of $\mathfrak{F}$ there exists $\sigma_s\in\sigma(\mathfrak{F})$, such that $\text{\rm Ind}_\sigma(k(\sigma_s))\le n-1$.
Let $h$ be a formation $\sigma$-function such that $h(\sigma_j)=f(\sigma_j)\cap \mathfrak{F}$ for all $\sigma_j\in\sigma$. Then $\mathfrak{F}=LF_{\sigma}(h)$ by Lemma~\ref{l3}. Since $\text{\rm Ind}_\sigma(\mathfrak{F})=n$, in view of Lemma~\ref{l5} $\text{\rm Ind}_\sigma(f(\sigma_j)\cap \mathfrak{F})=\text{\rm Ind}_\sigma(h(\sigma_j))\ge n-1$ for all $\sigma_j\in \sigma(\mathfrak{F})$. 
Note also that $F(\sigma_j)=\mathfrak{G}_{\sigma_j}h(\sigma_j)$ for all $\sigma_j\in\sigma$ by Lemma~\ref{l6}~(5).
In view of Lemma~\ref{l9} the formation $\mathfrak{G}_{\sigma_j}$ is $(n-1)$-multiply $\sigma$-local and also the product $\mathfrak{G}_{\sigma_j}h(\sigma_j)=F(\sigma_j)$ is an $(n-1)$-multiply $\sigma$-local formation for all $\sigma_j\in \sigma(\mathfrak{F})$.
Since $\text{\rm Ind}_\sigma(\mathfrak{F})=n$, there exists  $\sigma_j\in \sigma(\mathfrak{F})$ such that $\text{\rm Ind}_\sigma(F(\sigma_j))= n-1$.

Adequacy. Now suppose there exists $\sigma_i\in\sigma$ such that $\text{\rm Ind}_\sigma(F(\sigma_i))=n-1$, and for all $\sigma_j\in\sigma(\mathfrak{F})\setminus\{\sigma_i\}$ we have $\text{\rm Ind}_\sigma(F(\sigma_j))\ge n-1$. Suppose $\text{\rm Ind}_\sigma(\mathfrak{F})>n$. Then there exists a $\sigma$-local definition $h$ of $\mathfrak{F}$ such that $h(\sigma_i)$ is $n$-multiply $\sigma$-local. Then $h(\sigma_i)\cap\mathfrak F$ is $n$-multiply $\sigma$-local by Lemma~\ref{l5}. Now applying Lemmas~\ref{l3} and~ \ref{l6}~(5) we have $F(\sigma_i)=\mathfrak{G}_{\sigma_i}(h(\sigma_i)\cap\mathfrak F)$. In view of Lemma~\ref{l9} the formation $F(\sigma_i)=\mathfrak{G}_{\sigma_i}(h(\sigma_i)\cap\mathfrak F)$ is $n$-multiply $\sigma$-local. The resulting contradiction shows that $\text{\rm Ind}_\sigma(\mathfrak{F})=n$. The proposition is proved.

In the classical case, when $\sigma=\sigma^1$ from proposition~\ref{p1} we obtain

\begin{corollary}{\rm \cite[p. 27]{sk1}} \label{cp1}  
Let $f$ be the canonical screen of $\mathfrak{F}$.
Then $\text{\rm Ind}_l(\mathfrak{F})=n$ if and only if there is  $p\in\mathbb{P}$ such that $\text{\rm Ind}_l(f(p))=n-1$, and $\text{\rm Ind}_l(f(q))\ge n-1$ for all $q\in\pi(\mathfrak{F})\setminus\{p\}$.
\end{corollary}

Recall \cite[p.~28]{dh} that for any nonempty sets of primes $\pi$ and $\nu$ the symbol $O_{\pi,\nu} (G)$ denotes the characteristic subgroup of the group $G$, defined by the relation $O_{\pi,\nu} (G)/O_{\pi}(G)=O_{\nu} (G/O_{\pi}(G))$. 

We need the following special case of Lemma~14 from \cite{schsaf}.

\begin{lemma} {\rm \cite{schsaf}} \label{l17} 
Let $A$ be a non-identity $p$-group for some prime  $p$, $\pi$ and $\nu$ be nonempty sets of primes such that $p\in \nu$, $\pi\subseteq p'$. Then if \  $W_1=K\cdot B_1\le W=A\wr B=K\rtimes B$, where $K$ is the base group of the regular wreath product $W$, $B\ne 1$ and $B_1$ is a subgroup of $B$, then 
$$
O_\pi(W_1)=1, \ \ 
O_\nu (W_1)=K\cdot O_\nu(B_1),
$$
$$
O_{\nu ,\pi}(W_1)=K\cdot O_{\nu ,\pi}(B_1), \ \ 
O_{\pi,\nu} (W_1)=O_\nu (W_1)=K\cdot O_{\nu}(B_1).
$$
\end{lemma}

Recall that a formation $\sigma$-function $f$ is called {\sl $l^{\tau}_{\sigma_{n-1}}$-valued}, if $f(\sigma_{i})$ is a $\tau$-closed $(n-1)$-multiply $\sigma$-local formation for all $\sigma_i\in {\rm Supp}(f)$. Following \cite[p. 14]{sk1},  we denote by the symbol $(l^{\tau}_{\sigma_{n-1}})^\sigma$ a collection of all $\sigma$-local formations which have at least one $l^{\tau}_{\sigma_{n-1}}$-valued $\sigma$-local definition.  
 We will assume that by definition $ (1) \in(l^{\tau}_{\sigma_{n-1}})^\sigma $. 
 
\begin{theorem} \label{t2} 
$
(l_{\sigma_{n-1}}^\tau)^{\sigma} = l_{\sigma_n}^\tau
$
for any positive integer $n$.
\end{theorem}

Proof. Let us show first that $(l_{\sigma_{n-1}}^\tau)^{\sigma} \subseteq l_{\sigma_n}^\tau$. Let $\mathfrak{F}\in (l_{\sigma_{n-1}}^\tau)^{\sigma}$. Then $\mathfrak{F}=LF_{\sigma}(f)$, where $f(\sigma_i)\in l_{\sigma_{n-1}}^\tau$ for all $\sigma_i\in \text{Supp}(f)$. The latter means that $\mathfrak{F}$ is $n$-multiply $\sigma$-local by definition. 

Now we prove that $\mathfrak{F}$ is a $\tau$-closed formation. Let $G\in\mathfrak{F}$ and $H\in\tau(G)$. It is necessary to show that $H\in\mathfrak F$. Let $\sigma_i\in\sigma(H)$. Since $G\in \mathfrak F$, we have $G/ F_{\{\sigma_i\}}(G)\in f(\sigma_i)$. But the formation $f(\sigma_i)$ is $\tau$-closed and $H\in\tau(G)$, hence $H F_{\{\sigma_i\}}(G)/ F_{\{\sigma_i\}}(G)\in\tau(G/ F_{\{\sigma_i\}}(G))\subseteq f(\sigma_i).$

But then $H/ F_{\{\sigma_i\}}(G)\cap H\simeq H F_{\{\sigma_i\}}(G)/ F_{\{\sigma_i\}}(G)\in f(\sigma_i)$. Since $F_{\{\sigma_i\}}(G)\cap H\subseteq F_{\{\sigma_i\}}(H)$ we have $H/ F_{\{\sigma_i\}}(H)\in f(\sigma_i)$.
Hence, $H/ F_{\{\sigma_i\}}(H)\in f(\sigma_i)$ for all $\sigma_i\in\sigma(H)$. The latter implies $H\in\mathfrak{F}$. Hence, the formation $\mathfrak{F}$ is $\tau$-closed. Thus, $\mathfrak{F}\in l_{\sigma_n}^\tau$. Therefore $(l_{\sigma_{n-1}}^\tau)^{\sigma}\subseteq l_{\sigma_n}^\tau$.

Now let $\mathfrak{F}\in l_{\sigma_n}^\tau$ and $F$ be the canonical $\sigma$-local definition of $\mathfrak{F}$. Let us show that  $F(\sigma_i)\in l_{\sigma_{n-1}}^\tau$ for all $\sigma_i\in\sigma(\mathfrak{F})$. In view of the Proposition~\ref{p1} we have  $F(\sigma_i)\in l^{\sigma}_{n-1}$. Now we prove that   $F(\sigma_i)$ is $\tau$-closed. Let $G\in F(\sigma_i)$ and $H\in \tau(G)$. By induction on $\mid G\mid$, we show that $H\in F(\sigma_i)$. 

Suppose that $G$ has two different minimal normal subgroups $K_1$ and $K_2$.
Then $HK_j/K_j\in\tau(G/K_j)$, $j=1, 2$.
Hence, by induction $H/K_j\cap H\simeq HK_j/K_j\in F(\sigma_i)$, $j=1, 2$.
But then $H\simeq H/1 = H/(K_1\cap H)\cap (K_2\cap H)\in F(\sigma_i).$

Now let $K$ be the unique minimal normal subgroup of $G$. Suppose $K$ is a $\sigma_i$-group. Then $H^{F(\sigma_i)}\subseteq K\cap H\in\mathfrak{G}_{\sigma_i}$ since $H/K\cap H\simeq HK/K\in F(\sigma_i)$. Hence, $H\in \mathfrak{G}_{\sigma_i}F(\sigma_i)$. But $\mathfrak{G}_{\sigma_i}F(\sigma_i)=F(\sigma_i)$ in view of Lemma~\ref{l6}~(5). Therefore, $H\in F(\sigma_i)$.

Let $\mbox{O}_{\sigma_i}(G)=1$, $p\in\sigma_i$ and $V$ be a group of order $p$. We denote by  $W$ the regular wreath product $V\wr G$ of the groups $V$ and $G$. Then $W=P\rtimes G$, where $P$ is the base group of $W$. Since $\mbox{O}_{\sigma_i}(G)=1$ and $P\in\mathfrak G_{\sigma_i}$ we have $\mbox{O}_{\sigma_i}(W)=\mbox{O}_{\sigma_i}(W)\cap PG= P(\mbox{O}_{\sigma_i}(W)\cap G)=P$. But then $W/P\simeq G\in F(\sigma_i)\subseteq\mathfrak F$. Hence, $W\in\mathfrak{F}$ by Lemma~\ref{l1}. Now if $\varphi : W \to G$ is a natural epimorphism of $W$ on $G$, then, obviously, $H^{\varphi^{-1}} = PH$. Hence, $PH\in\tau(W)$ and $PH\in\mathfrak{F}$. Therefore, $PH/ F_{\{\sigma_i\}}(PH)\in F(\sigma_i)$.

In view of Lemma~\ref{l17} (for $\nu = \sigma_i$, $\pi=\sigma_i'$\,) we have
$$
\mbox{O}_{\sigma_i' ,\,\sigma_i}(PH)=F_{\{\sigma_i\}}(PH) = \mbox{O}_{\sigma_i}(PH)= P\mbox{O}_{\sigma_i}(H).
$$
Consequently,
$$
PH/ F_{\{\sigma_i\}}(PH) = PH/P\mbox{O}_{\sigma_i}(H)\simeq H/\mbox{O}_{\sigma_i}(H)(P\cap H) = H/\mbox{O}_{\sigma_i}(H)\in F(\sigma_i).
$$
But then $H\in F(\sigma_i)$ by Lemma~\ref{l1}.
Thus, $F(\sigma_i)$ is a $\tau$-closed formation. Therefore, $F(\sigma_i)\in l_{\sigma_{n-1}}^\tau$.
Hence, $\mathfrak{F}\in (l_{\sigma_{n-1}}^\tau)^{\sigma}$ and  $l_{\sigma_n}^\tau\subseteq (l_{\sigma_{n-1}}^\tau)^{\sigma}$. The theorem is proved.

\begin{corollary} \label{c11t2} 
Let $\mathfrak{F}=LF_{\sigma}(F)$ be a $\tau$-closed $n$-multiply $\sigma$-local formation, where $n\ge 1$ and $F$ is the canonical $\sigma$-local definition of $\mathfrak F$. Then $F(\sigma_i)\in l^{\tau}_{\sigma_{n-1}}$ for all $\sigma_i\in\sigma(\mathfrak F).$
\end{corollary}

If $\tau$ is a trivial subgroup functor from Theorem~\ref{t2} we obtain

\begin{corollary}\label{c12t2} 
$
(l^{\sigma}_{n-1})^{\sigma} = l^{\sigma}_n
$ 
for any positive integer $n$.
\end{corollary}

\begin{corollary}\label{c13t2} 
Let $\mathfrak{F}=LF_{\sigma}(F)$ be an $n$-multiply $\sigma$-local formation, where $n\ge 1$ and  $F$ is the canonical $\sigma$-local definition of $\mathfrak F$. Then $F(\sigma_i)\in l^{\sigma}_{n-1}$ for all $\sigma_i\in\sigma(\mathfrak F).$
\end{corollary}

In the classical case, when $\sigma=\sigma^1$ from Theorem~\ref{t2} we obtain the following well-known results

\begin{corollary} {\rm \cite[p. 29]{sk1}}\label{c14t2} 
$
(l^\tau_{n-1})^l = l^\tau _n
$
for any positive integer $n$.
\end{corollary}

\begin{theorem} \label{t2a} 
 Let $\mathfrak{F}$ be an $n$-multiply $\sigma$-local formation $(n\ge 1)$ and let $F$ be the canonical $\sigma$-local definition of $\mathfrak F$. Then $\mathfrak{F}$ is $\tau$-closed if and only if 
 $F(\sigma_i)$ is a $\tau$-closed formation for all $\sigma_i\in\sigma(\mathfrak F)$.
\end{theorem}
 
Proof. Let $\mathfrak{F}$ be a $\tau$-closed $n$-multiply $\sigma$-local formation, where $n\ge 1$. Then $F(\sigma_i)$ is a $\tau$-closed formation for all $\sigma_i\in\sigma(\mathfrak F)$ by Corollary~\ref{c11t2}. 
On the other hand, if $F(\sigma_i)$ is a $\tau$-closed formation for all $\sigma_i\in\sigma(\mathfrak F)$, then in view of Lemma~\ref{l16a} $\mathfrak{F}$ is a $\tau$-closed formation as well. The theorem is  proved.
 
Let $\tau$ be a single subgroup functor, that is, $\tau(G)$ is the collection of all subgroups of $G$ for any group $G$. Then from Theorem~\ref{t2a} we obtain the following results 

\begin{corollary} \label{c14at2a}  
Let $\mathfrak{F}$ be an $n$-multiply $\sigma$-local formation $(n\ge 1)$ and let $F$ be the canonical $\sigma$-local definition of $\mathfrak F$. Then $\mathfrak{F}$ is $\tau$-closed if and only if 
 $F(\sigma_i)$ is a $\tau$-closed formation for all $\sigma_i\in\sigma(\mathfrak F)$.
 \end{corollary}
 
\begin{corollary} {\rm \cite[p. 41]{sh}}\label{c14bt2a} 
Let $f$ be the maximal integrated local screen of the formation  $\mathfrak F$. The formation $\mathfrak{F}$ is hereditary if and only if for any prime $p$ the formation $f(p)$ is hereditary. 
\end{corollary}  

Let $\tau(G)=\text{\rm S}_n(G)$ be the collection of all normal subgroups of $G$ for any group $G$ \cite[p.17]{sk1}. Then from Theorem~\ref{t2a} we obtain the following results

\begin{corollary} \label{c14ct2a} 
 Let $\mathfrak{F}$ be an $n$-multiply $\sigma$-local formation 
 $(n\ge 1)$ and let $F$ be the canonical $\sigma$-local definition of $\mathfrak F$. Then $\mathfrak{F}$ is normal hereditary if and only if 
 $F(\sigma_i)$ is a normal hereditary formation for all 
 $\sigma_i\in\sigma(\mathfrak F)$.
\end{corollary}

\begin{corollary} {\rm \cite[p. 41]{sh}} \label{c14dt2a} 
Let $f$ be the maximal integrated local screen of the formation  $\mathfrak F$. The formation $\mathfrak{F}$ is normal hereditary if and only if for any prime $p$ the formation $f(p)$ is normal hereditary. 
\end{corollary}

\begin{lemma} \label{l18} 
Let $\mathfrak{F}=LF_\sigma (f)$ and $G$ be a group with $G^{\mathfrak{F}}\ne 1$. Then there exists $\sigma_i\in\sigma(G^{\mathfrak{F}})$, such that $G/ F_{\{\sigma_i\}}(G)\notin f(\sigma_i)$.
\end{lemma}

Proof. Since $G\notin\mathfrak F$, there exists $\sigma_i\in\sigma(G)$, such that $G/ F_{\{\sigma_i\}}(G)\notin f(\sigma_i)$. 
Note that if $\sigma_j\in\sigma(G)\setminus\sigma(G^{\mathfrak{F}})$, then $G^{\mathfrak F}\subseteq \mbox{O}_{\sigma_j'}(G)\subseteq F_{\{\sigma_j\}}(G)$. Therefore, $F_{\{\sigma_j\}}(G/G^{\mathfrak{F}})= F_{\{\sigma_j\}}(G)/G^{\mathfrak{F}}$. 
Since $G/G^{\mathfrak{F}}\in\mathfrak{F}$, $(G/G^{\mathfrak{F}})/ F_{\{\sigma_j\}}(G/G^{\mathfrak{F}})\in f(\sigma_j)$ for all $\sigma_j\in\sigma(G)\setminus\sigma(G^{\mathfrak{F}})$. Hence, 
$$
(G/G^{\mathfrak{F}})/ F_{\{\sigma_j\}}(G/G^{\mathfrak{F}})=(G/G^{\mathfrak{F}})/( F_{\{\sigma_j\}}(G)/G^{\mathfrak{F}})\simeq G/ F_{\{\sigma_j\}}(G) \in  f(\sigma_j).
$$
Therefore, $\sigma_i\in\sigma(G^{\mathfrak{F}})$. The lemma is proved. 
 
\begin{theorem} \label{t3} 
Let $\mathfrak{F}$ be a nonempty $\tau$-closed formation.
Then $\mathfrak{F}$ is an $n$-multiply $\sigma$-local $(n\ge 1)$ if and only if  $\mathfrak{G}_{\sigma_i}\mathfrak{F}_{\sigma_{n-1}}^{\tau}(\sigma_i)\subseteq\mathfrak{F}$ for all $\sigma_i\in\sigma$.
\end{theorem}

Proof. Let $\mathfrak{F}$ be $n$-multiply $\sigma$-local, $F$ be the canonical $\sigma$-local definition of $\mathfrak{F}$.
Since $\mathfrak{F}\in l_{\sigma_n}^\tau$,  a formation $F(\sigma_i)$ is $\tau$-closed $(n-1)$-multiply $\sigma$-local by 
Corollary~\ref{c11t2}. Hence,
$$
\mathfrak{F}_{\sigma_{n-1}}^\tau (\sigma_i) = l_{\sigma_{n-1}}^\tau\mbox{form\,}(\mathfrak{F}(\sigma_i)) \subseteq F(\sigma_i).
$$
But $\mathfrak{G}_{\sigma_i}F(\sigma_i) = F(\sigma_i)$ by Lemma~\ref{l2}. Consequently,
$$
\mathfrak{G}_{\sigma_i}\mathfrak{F}_{\sigma_{n-1}}^\tau (\sigma_i)\subseteq \mathfrak{G}_{\sigma_i}F(\sigma_i) = F(\sigma_i) \subseteq \mathfrak{F}.
$$

Now suppose that $\mathfrak{G}_{\sigma_i}\mathfrak{F}_{\sigma_{n-1}}^\tau (\sigma_i)\subseteq \mathfrak{F}$ for all $\sigma_i\in\sigma$. Let $f$ be a formation $\sigma$-function such that $f(\sigma_i) = \mathfrak{G}_{\sigma_i}\mathfrak{F}_{\sigma_{n-1}}^\tau (\sigma_i)$ for all $\sigma_i\in\sigma$. Let us show that $\mathfrak{F}=LF_{\sigma}(f)$. Suppose that $\mathfrak{F}\not\subseteq LF_{\sigma}(f)$ and $G$ is a group of minimal order in $\mathfrak{F}\setminus LF_{\sigma}(f)$. Then $G$ is a monolithic group with monolith $R =G^{LF_{\sigma}(f)}$. Let $\sigma_i\in\sigma(R)$. If $R$ is not a $\sigma$-primary group, then $F_{\{\sigma_i\}}(G)=1$. Consequently,
$$
G\simeq G/F_{\{\sigma_i\}}(G)\in\mathfrak{F}(\sigma_i)\subseteq \mathfrak{G}_{\sigma_i}\mathfrak{F}_{\sigma_{n-1}}^\tau (\sigma_i)=f(\sigma_i).
$$
Hence, $G/F_{\{\sigma_i\}}(G)\in f(\sigma_i)$ for all $\sigma_i\in\sigma(R)$. The latter contradicts Lemma~\ref{l18}. Therefore, $R$ is a $\sigma_i$-group. Then $\mbox{O}_{\sigma_i'}(G)=1$ and $F_{\{\sigma_i\}}(G)=\mbox{O}_{\sigma_i}(G)$. 
Hence,
$$
G/\mbox{O}_{\sigma_i}(G)=G/F_{\{\sigma_i\}}(G)\in\mathfrak{F}(\sigma_i)\subseteq \mathfrak{G}_{\sigma_i}\mathfrak{F}_{\sigma_{n-1}}^\tau (\sigma_i)=f(\sigma_i).
$$
But then, $G\in LF_{\sigma}(f)$ in view of Lemma~\ref{l1}. A contradiction. Therefore, $\mathfrak{F}\subseteq LF_{\sigma}(f)$.

Now suppose that $LF_{\sigma}(f)\not\subseteq\mathfrak{F}$ and let $G$ be a group of minimal order in $LF_{\sigma}(f)\setminus\mathfrak{F}$. Then $G$ is monolithic and its monolith $R=G^{\mathfrak{F}}$. Let $\sigma_i\in\sigma(R)$. Suppose that $R$ is not a $\sigma$-primary group. Then $F_{\{\sigma_i\}}(G)=1$. Hence, since $G\in LF_{\sigma}(f)$ we have
$$
G\simeq G/F_{\{\sigma_i\}}(G)\in f(\sigma_i)=\mathfrak{G}_{\sigma_i}\mathfrak{F}_{\sigma_{n-1}}^\tau (\sigma_i)\subseteq\mathfrak{F}.
$$
This contradicts the choice of $G$. Therefore, $R$ is a $\sigma_i$-group. Then $\mbox{O}_{\sigma_i}(G)=F_{\{\sigma_i\}}(G)$, and so
$$
G/\mbox{O}_{\sigma_i}(G)=G/F_{\{\sigma_i\}}(G)\in f(\sigma_i)=\mathfrak{G}_{\sigma_i}\mathfrak{F}_{\sigma_{n-1}}^\tau (\sigma_i).
$$ Therefore,
$$
G\in\mathfrak{G}_{\sigma_i}\mathfrak{F}_{\sigma_{n-1}}^\tau (\sigma_i)\subseteq\mathfrak{F}.
$$
A contradiction. Therefore $LF_{\sigma}(f)\subseteq\mathfrak{F}$. Thus, $LF_{\sigma}(f)=\mathfrak{F}$.
The theorem is proved.

\begin{corollary} \label{c15t3} 
Let $\mathfrak{F}$ be a nonempty $\tau$-closed formation.
Then $\mathfrak{F}$ is $\sigma$-local if and only if 
$\mathfrak{G}_{\sigma_i}\mathfrak{F}^{\tau}(\sigma_i)\subseteq\mathfrak{F}$ for all $\sigma_i\in\sigma$.
\end{corollary}

If $\tau$ is a trivial subgroup functor, then from Theorem~\ref{t3} we obtain

\begin{corollary} \label{c16t3} 
Let $\mathfrak{F}$ be a nonempty formation.
Then $\mathfrak{F}$ is $n$-multiply $\sigma$-local $(n\ge 1)$ if and only if  
$\mathfrak{G}_{\sigma_i}\mathfrak{F}^{\sigma}_{n-1}(\sigma_i)\subseteq\mathfrak{F}$ for all $\sigma_i\in\sigma$.
\end{corollary}

In the classical case, when $\sigma=\sigma^1$, as corollaries of Theorem~\ref{t3} we obtain the following well-known result

\begin{corollary} {\rm \cite[p. 31]{sk1}} \label{c17t3} 
Let $\mathfrak{F}$ be a nonempty  $\tau$-closed formation.
Then $\mathfrak{F}$ $n$-multiply local $(n\ge 1)$ if and only if 
$\mathfrak{N}_p\mathfrak{F}_{n-1}^{\tau}(p)\subseteq\mathfrak{F}$  for all $p\in\mathbb{P}$.
\end{corollary}

\begin{theorem} \label{t4} 
Let $\mathfrak{X}$ be some nonempty collection of groups, $n\ge 1$. Then
$$
l_{\sigma_n}^\tau\mbox{\rm form\,}\mathfrak{X}=\mbox{\rm form\,}(\cup_{\sigma_i\in\sigma}\mathfrak{G}_{\sigma_i}\mathfrak{X}_{\sigma_{n-1}}^\tau(\sigma_i)).
$$
\end{theorem}

Proof. Let $ \mathfrak{F}=l_{\sigma_n}^\tau\mbox{\rm form\,}\mathfrak{X} $ and  $\mathfrak{M}=\mbox{\rm form\,}(\cup_{\sigma_i\in\sigma}\mathfrak{G}_{\sigma_i}\mathfrak{X}_{\sigma_{n-1}}^\tau(\sigma_i))$. In view of Theorem~\ref{t3}  we have $\mathfrak{G}_{\sigma_i}\mathfrak{F}_{\sigma_{n-1}}^{\tau}(\sigma_i)\subseteq\mathfrak{F}$ for all $\sigma_i\in\sigma$. Hence, $\mathfrak{M}\subseteq\mathfrak{F}$.
Moreover, if $A\in\mathfrak{X}$, then
$$
A/\mbox{O}_{\sigma_i', \sigma_i}(A)=A/F_{\{\sigma_i\}}(A)\in\mathfrak{X}(\sigma_i)\subseteq\mathfrak{X}_{\sigma_{n-1}}^\tau(\sigma_i).
$$
Therefore,
$$
A/\mbox{O}_{\sigma_i'}(A)\in\mathfrak{G}_{\sigma_i}\mathfrak{X}_{\sigma_{n-1}}^\tau(\sigma_i)\subseteq\mathfrak{M}.
$$
Since in this case $\mathfrak M$ is a formation and $\cap_{\sigma_i\in\sigma(A)}\mbox{O}_{\sigma_i'}(A)=1$, we have
$$
A\simeq A/(\cap_{\sigma_i\in\sigma(A)}\mbox{O}_{\sigma_i'}(A))\in\mathfrak{M}.
$$

Hence, $\mathfrak{X}\subseteq\mathfrak{M}$. Thus, there are inclusions $\mathfrak{X}\subseteq\mathfrak{M}\subseteq\mathfrak{F}$. Since $\mathfrak{F}=l_{\sigma_n}^\tau\mbox{\rm form\,}\mathfrak{X},$ it is sufficient to establish that the formation $\mathfrak{M}$ is $\tau$-closed $n$-multiply $\sigma$-local.

Let us show that $\mathfrak{M}$ is $\tau$-closed. Indeed, since for any $\sigma_i\in\sigma$ the formation $\mathfrak{G}_{\sigma_i}$ is hereditary, and the formation $\mathfrak{X}_{\sigma_{n-1}}^\tau(\sigma_i)$ is a $\tau$-closed, we have $\mathfrak{G}_{\sigma_i}\mathfrak{X}_{\sigma_{n-1}}^\tau(\sigma_i)$ is a $\tau$-closed formation by Lemma~\ref{l8}. Hence, in view of  Lemma~\ref{l16} $\mathfrak{M}$ is $\tau$-closed.

Now we prove that $\mathfrak{M}$ is an $n$-multiply $\sigma$-local formation.
In view of Theorem~\ref{t3}, it suffices to show that $\mathfrak{G}_{\sigma_i}\mathfrak{M}_{\sigma_{n-1}}^\tau(\sigma_i)\subseteq\mathfrak{M}$ for all $\sigma_i\in\sigma$. 
Since $\mathfrak{F}=l_{\sigma_n}^\tau\mbox{\rm form\,}\mathfrak{X},$   we have $\mathfrak{X}_{\sigma_{n-1}}^\tau(\sigma_i)=\mathfrak{F}_{\sigma_{n-1}}^\tau(\sigma_i)$ for all $\sigma_i\in\sigma$ by Theorem~\ref{t1}. 
On the other hand, $\mathfrak{M}\subseteq\mathfrak{F}$, then, obviously, $\mathfrak{M}_{\sigma_{n-1}}^\tau(\sigma_i)\subseteq \mathfrak{F}_{\sigma_{n-1}}^\tau(\sigma_i)$. Therefore, $\mathfrak{M}_{\sigma_{n-1}}^\tau(\sigma_i)\subseteq \mathfrak{X}_{\sigma_{n-1}}^\tau(\sigma_i)$. Hence,  
$$
\mathfrak{G}_{\sigma_i}\mathfrak{M}_{\sigma_{n-1}}^\tau(\sigma_i)\subseteq \mathfrak{G}_{\sigma_i}\mathfrak{X}_{\sigma_{n-1}}^\tau(\sigma_i)\subseteq \mathfrak M.
$$
Thus, $\mathfrak M$ is $n$-multiply $\sigma$-local. The theorem is proved.

\begin{corollary} \label{c18t4} 
Let $\mathfrak{X}$ be nonempty collection of groups. Then
$
l_{\sigma}^\tau\mbox{\rm form\,}\mathfrak{X}=\mbox{\rm form\,}(\cup_{\sigma_i\in\sigma}\mathfrak{G}_{\sigma_i}\mathfrak{X}^\tau(\sigma_i)).
$
\end{corollary}

In the case when $\tau$ is a trivial subgroup functor we have 

\begin{corollary} \label{c19t4} 
Let $\mathfrak{X}$ be nonempty collection of groups, $n\ge 1$. Then
$$
l^{\sigma}_n\mbox{\rm form\,}\mathfrak{X}=\mbox{\rm form\,}(\cup_{\sigma_i\in\sigma}\mathfrak{G}_{\sigma_i}\mathfrak{X}^{\sigma}_{n-1}(\sigma_i)).
$$
\end{corollary}
 
In the classical case from Theorem~\ref{t4} we obtain the following well-known result
 
\begin{corollary} {\rm \cite[p. 33]{sk1}} \label{c20t4} 
Let $\mathfrak{X}$ be some nonempty collection of groups, $n\ge 1$. Then
$$
l_{n}^\tau\mbox{\rm form\,}\mathfrak{X}=\mbox{\rm form\,}(\cup_{p\in\mathbb{P}}\mathfrak{N}_{p}\mathfrak{X}_{n-1}^\tau(p)).
$$
\end{corollary}

\section{\bf Lattice of all $\tau$-closed $n$-multiply $\sigma$-local formations}

Recall that a nonempty  collection of formations $\theta$ is called a complete lattice of formations~\cite[p. 12]{sk1},  if the intersection of any set of formations from $\theta$ again belongs to $\theta$ and the set $\theta$ contains a formation $\mathfrak{F}$, such that $\mathfrak{H}\subseteq\mathfrak{F}$ for all $\mathfrak{H}\in\theta$. 
                       
\begin{lemma} \label{l19} 
If the class of groups $\mathfrak{F}_{j}$ is a $\tau$-closed $n$-multiply $\sigma$-local formation for all $j\in J$, then the class $\mathfrak{F}=\cap _{j\in J}\mathfrak{F}_{j}$ is also a $\tau$-closed $n$-multiply $\sigma$-local formation.
\end{lemma} 

Proof. By Lemma~\ref{l5} the formation $\mathfrak F$ is $n$-multiply $\sigma$-local. Therefore, it suffices to show that $\cap _{j\in J}\mathfrak{F}_{j}$ is a $\tau$-closed formation. Indeed, if $G\in\cap _{j\in J}\mathfrak{F}_{j}$ and $H\in\tau(G)$, then $G\in\mathfrak{F}_{j}$ for any $j\in J$. Then, since the formation $\mathfrak F_j$ is $\tau$-closed, we have $H\in\mathfrak F_j $. Hence, $H\in\cap _{j\in J}\mathfrak{F}_{j}=\mathfrak F$. Thus, $\mathfrak{F}$ is a $\tau$-closed $n$-multiply $\sigma$-local formation. The lemma is proved.

Recall that if $\{ \mathfrak{F}_j \mid j\in J \}$ is some collection of $\tau$-closed $n$-multiply $\sigma$-local formations, then we put
$$
\vee ^{\tau}_{\sigma_n}(\mathfrak{F}_j \mid j\in J )=l^{\tau}_{\sigma_n} {\rm form\,}( \cup_ {j\in J}\mathfrak{F}_j).
$$
In particular, if $\tau$ is a trivial subgroup functor,  then 
$$
\vee ^{\sigma}_n(\mathfrak{F}_j \mid j\in J )=l^{\sigma}_n {\rm form\,}( \cup_ {j\in J}\mathfrak{F}_j).
$$
 
\begin{theorem} \label{t5}  
The collection $l^{\tau}_{\sigma_n}$ of all $\tau$-closed $n$-multiply $\sigma$-local formations is partially ordered with respect to the inclusion. Moreover, $l^{\tau}_{\sigma_n}$ forms a complete lattice of formations in which, for any set $\{ \mathfrak{F}_j \mid j\in J \} \subseteq l^{\tau}_{\sigma_n}$, the intersection $\cap _{j\in J}\mathfrak{F}_j$ is the greatest lower bound and $\vee ^{\tau}_{\sigma_n}(\mathfrak{F}_j \mid j\in J )$ is the smallest upper bound $\{ \mathfrak{F}_j \mid j\in J \}$ in $l^{\tau}_{\sigma_n}$.
\end{theorem}

Proof. The first statement is obvious. Note also that the class $\mathfrak{G}$ of all groups is, obviously, a $\tau$-closed formation for any subgroup functor $\tau$. In addition, by Lemma~\ref{l7} the formation $\mathfrak{G}$ is totally $\sigma$-local, and hence, $\mathfrak{G}$ is $n$-multiply $\sigma$-local for every non-negative integer $n$. Therefore, $\mathfrak{G}\in l^{\tau}_{\sigma_n}$. So it is the largest element in $l^{\tau}_{\sigma_n}$. On the other hand, $\cap _{j\in J}{\mathfrak{F}}_{j}\in l^{\tau}_{\sigma_n}$ by Lemma~\ref{l19}. Therefore, $\cap _{j\in J}{\mathfrak{F}}_{j}$ is the greatest lower bound $\{ {\mathfrak{F}}_{j} \mid j\in J \}$ in $ l^{\tau}_{\sigma_n}$. The latter also means that $\vee ^{\tau}_{\sigma_n}(\mathfrak{F}_j \mid j\in J )$ is the smallest upper bound $\{ {\mathfrak{F}}_{j} \mid j\in J \}$ in $ l^{\tau}_{\sigma_n}$.The theorem is proved.

In the case when $\tau$ is a trivial subgroup functor from Theorem~\ref{t5} we obtain

\begin{corollary} {\rm \cite{chsafsk}} \label{c21t5} 
The collection $l^{\sigma}_n$ of all $n$-multiply $\sigma$-local formations is partially ordered with respect to the inclusion. Moreover, $l^{\sigma}_n$  is a complete lattice in which, for any set  $\{ \mathfrak{F}_j \mid j\in J \} \subseteq l^{\sigma}_n$, the intersection $\cap _{j\in J}\mathfrak{F}_j$ is the greatest lower bound and $\vee ^{\sigma}_n(\mathfrak{F}_j \mid j\in J )$ is the smallest upper bound $\{ \mathfrak{F}_j \mid j\in J \}$ in $l^{\sigma}_n$.
\end{corollary}

\begin{lemma} \label{l21} 
Let $f_j$ be the smallest $l^{\tau}_{\sigma_{n-1}}$-valued $\sigma$-local definition of $\mathfrak{F}_j$, $j\in J$. Then $\vee^{\tau}_{\sigma_{n-1}}(f_j\mid j\in J)$ is the smallest $l^{\tau}_{\sigma_{n-1}}$-valued  $\sigma$-local definition of $\mathfrak{F}=\vee^{\tau}_{\sigma_{n}}(\mathfrak{F}_j\mid j\in J)$.
\end{lemma} 

Proof. Let $m$ be the smallest $l^{\tau}_{\sigma_{n-1}}$-valued definition of $\mathfrak F,$ $f=\vee^{\tau}_{\sigma_{n-1}}(f_{j} \mid j\in J)$ and $\Pi =\sigma (\cup _{j\in J} {\rm {\mathfrak F}}_{j} )=\cup _{j\in J} \sigma ({\rm {\mathfrak F}}_{j} ).$ Then $\sigma ({\rm {\mathfrak F}})=\Pi $ by Theorem~\ref{t1}. 
Now show that $m(\sigma_i )=f(\sigma_i )$ for all $\sigma_i \in \sigma.$

Let $\sigma_i \in \sigma \backslash \Pi .$ Then for any $j\in J$ we have $f_{j} (\sigma_i )=\emptyset $ by Theorem~\ref{t1}. 
Hence, $f(\sigma_i )=\emptyset .$ On the other hand, since $\sigma_i\notin\sigma(\mathfrak F)$, by  Theorem~\ref{t1} $m(\sigma_i )=\emptyset.$ So, $m(\sigma_i)=f(\sigma_i ).$

Now suppose $\sigma_i \in \Pi .$ Then there exists $j_{t} \in J$ such that $\sigma_i \in \sigma ({\rm {\mathfrak F}}_{j_{t} } ).$ From Theorem~\ref{t1} it follows that $f_{j_{t} } (\sigma_i )\ne \emptyset $ and 
$$
m(\sigma_i )=l^{\tau}_{\sigma_{n-1}}{\rm form\,}(G/F_{\{ \sigma_i \} } (G)\mid G\in \cup _{j\in J} {\rm {\mathfrak F}}_{j} )=$$
$$=l^{\tau}_{\sigma_{n-1}}{\rm form\,}(\cup _{j\in J} l^{\tau}_{\sigma_{n-1}} {\rm form\,}(G/F_{\{ \sigma_i \} } (G)\mid G\in {\rm {\mathfrak F}}_{j} ))=$$ 
$$=l^{\tau}_{\sigma_{n-1}} {\rm form\,}(\cup _{j\in J} f_{j} (\sigma_i )\mid j\in J)=(\vee^{\tau}_{\sigma_{n-1}}  (f_{j} \mid j\in J))(\sigma_i )=f(\sigma_i ).$$ 
Therefore, $m(\sigma_i )=f(\sigma_i )$ for all $\sigma_i \in \sigma .$  The lemma is proved.

In particular, in the case when $\tau$ is a trivial subgroup functor, we have
  
\begin{lemma} \label{l22} 
Let $f_j$ be the smallest $l^{\sigma}_{n-1}$-valued $\sigma$-local definition of the formation $\mathfrak{F}_j$, $j\in J$. Then $\vee^{\sigma}_{n-1}(f_j\mid j\in J)$ is the smallest $l^{\sigma}_{n-1}$-valued $\sigma$-local definition of the formation $\mathfrak{F}=\vee^{\sigma}_{n}(\mathfrak{F}_j\mid j\in J)$.
\end{lemma} 

\begin{lemma} \label{l23}  
Let $A$ be a monolithic group such that ${\rm Soc}(A)$ is not a $\sigma$-primary group, and let $\mathfrak M$ be some $\tau$-closed semiformation. Then if $A\in l^{\tau}_{\sigma_n} \text{\rm form\,}{\mathfrak M}$, then $A\in \mathfrak M.$
\end{lemma}

Proof. Let $K=\mbox{Soc}(A)$ and $m$ be the smallest $l^{\tau}_{\sigma_{n-1}}$-valued $\sigma$-local definition of  $l^{\tau}_{\sigma_n} \text{\rm form\,}{\mathfrak M}$.
We prove this lemma by induction on $n$. If $n=0$, then $ A \in l_{\sigma_0}^\tau \mbox{form\,}\mathfrak{M}=\tau \mbox{form\,}\mathfrak{M}.$
Suppose $A \notin \mathfrak{M}$. Then according to Lemma~\ref{l13} in $\tau \mbox{form\,}\mathfrak{M}$ there is a group $H$ with such normal subgroups $N, M, N_1, \ldots, N_t, M_1, \ldots,M_t,$ $(t \ge 2)$, that the following statements are true:
1) $A \simeq H/N$ and $M/N=\mbox{Soc}(H/N)$;
2) $N_1\cap \ldots \cap N_t=1;$
3) $H/N_j$ is a monolithic $\mathfrak{M}$-group with a monolith $M_j/N_j$, which is $H$-isomorphic to $M/N$;
4) $M_1\cap \ldots \cap M_t\subseteq M$.

By hypothesis $K$ is not a $\sigma$-primary unique minimal normal subgroup of $A$, so $K$ is a non-abelian group and $\mbox{C}_G(K)=1$. Since in this case $A\simeq H/N$, $\mbox{C}_{H/N}(M/N)=N/N$. Therefore, $\mbox{C}_H(M/N)=N$. 
Hence, $N_j \subseteq N$. But then, $(H/N_j)/(N/N_j)\simeq H/N\in \mathfrak M$, since $H/N_j$ is a monolithic $\mathfrak{M}$-group. Therefore, $A \simeq H/N \in \mathfrak{M}$. The resulting contradiction completes the proof of the lemma for $n=0$.

Now assume that $n>0$ and this lemma is true for $n-1$. Let $\sigma_i \in \sigma(K)$. Then $F_{\{\sigma_i\}}(A)=1$. Since $A \in l^\tau_{\sigma_n} \mbox{form\,}\mathfrak{M}$, by Theorem~\ref{t1} we have
$$
A \simeq A/F_{\{\sigma_i\}}(A) \in m(\sigma_i)\subseteq l_{\sigma_{n-1}}^\tau \mbox{form\,}\mathfrak{M}.
$$
Hence, by induction, $A \in \mathfrak{M}$.
The lemma is proved.

\begin{lemma}\label{l24} 
Let $N_1\times\ldots\times N_t=\mbox{Soc}(G)$, where $t>1$ and $G$ be a group with $\mbox{O}_{\sigma_i}(G)=1$. Let $M_k$ be the largest normal subgroup in $G$ containing $N_1\times\ldots\times N_{k-1}\times N_{k+1}\times\ldots\times N_t$, but not containing $N_k$, $k=1,\ldots, t$. Then the following statements are true:

(1) for any $k\in\{1,\ldots ,t\}$ the group $G/M_k$ is monolithic and its monolith $N_kM_k/M_k$ is $G$-isomorphic to $N_k$ and $\mbox{O}_{\sigma_i}(G/M_k)=1$;

(2) $M_1\cap\ldots\cap M_t=1$.
\end{lemma}

Proof. Let us prove (1). Suppose that the group $G/M_k$ is not monolithic. Then $G/M_k$ contains a minimal normal subgroup $L/M_k$, different from $N_kM_k/M_k$. Hence, $N_k\not\subseteq L$ and $N_1\times\ldots\times N_{k-1}\times N_{k+1}\times\ldots\times N_t\subseteq M_k\subseteq L$.
But then, according to the definition of the subgroup $M_k$ we have $L\subseteq M_k$. Therefore, $L=M_k$, contradiction.
Thus, the group $G/M_k$ is monolithic and $N_kM_k/M_k=\mbox{Soc}(G/M_k)$. Since $M_k\cap N_k=1$, there is a $G$-isomorphism $N_k\simeq N_kM_k/M_k$. Therefore, in view of the condition $\mbox{O}_{\sigma_i}(G)=1$ we obtain that $\mbox{O}_{\sigma_i}(G/M_k)=1$. Thus, statement (1) is proved.

Now we prove (2). Suppose that $M_1\cap\ldots\cap M_t\ne 1$ and let $D$ be the minimal normal subgroup of the group $G$, contained in $M_1\cap\ldots\cap M_t$. Then, for all $k=1,\ldots, t,$ obviously $D\not= N_k$ holds.
Since $\cap_{i=1}^{t}(N_1\times\ldots\times N_{k-1}\times N_{k+1}\times\ldots\times N_t)=1$, there exists $k\in\{1,\ldots , t\}$, such that $D\not\subseteq B=N_1\times\ldots\times N_{k-1}\times N_{k+1}\times\ldots\times N_t$.
Hence, $DB=\mbox{Soc}(G)\subseteq M_k$. But then $N_k\subseteq M_k$. The latter contradicts the hypothesis of the lemma. So, $M_1\cap\ldots\cap M_t=1$. The lemma is proved.

\begin{lemma} \label{l25} 
Let $\mathfrak{M}$ be a $\tau$-closed semiformation and $A\in l^{\tau}_{\sigma_n}\mbox{\rm form\,}\mathfrak{M}$. Then if $\mbox{O}_{\sigma_i}(A)=1$, then $A\in l^{\tau}_{\sigma_n}\mbox{\rm form\,}\mathfrak{X}$, where $\mathfrak{X}=\{G/\mbox{O}_{\sigma_i}(G)\mid G\in\mathfrak{M}\}$.
\end{lemma}

Proof. It is clear that if $A\in\mathfrak{M}$, then the lemma is true. Let $A\notin\mathfrak{M}$.
 First suppose that $A$ is a monolithic group and let $R=Soc(A)$. Let $n=0$. Then $ A\in l_{\sigma_0}^{\tau}\mbox{form\,}\mathfrak{M}=\tau\mbox{form\,}\mathfrak{M}.$
Hence, in view of Lemma~\ref{l13}, in $\tau\mbox{form\,}\mathfrak{M}$ there exists a group $H$ with normal subgroups $N, M, N_1, \ldots, N_t, M_1, \ldots, M_t,$ $(t\ge 2)$, such that the following conditions are hold:
1)~$H/N\simeq A$, $M/N=\mbox{Soc}(H/N)$;
2)~$N_1\cap\ldots\cap N_t=1$;
3) $H/N_k$ is a monolithic $\mathfrak{M}$-group with a monolith $M_k/N_k$, which is $H$-isomorphic to $M/N$.
It follows from conditions~1) and~3) that $\mbox{O}_{\sigma_i}(H/N_k)=1$ since $\mbox{O}_{\sigma_i}(A)=1$. Therefore, $\mathfrak{M}$-groups $H/N_1,\ldots, H/N_t$ belong to $\mathfrak X$. In addition, it follows from conditions~2) and~3) that $H$ is a subdirect product of the groups $H/N_1,\ldots, H/N_t$. Hence, $H\in\mbox{R}_0\{H/N_1,\ldots , H/N_t\},$ so 
$$
A\simeq H/N\in\mbox{QR}_0\{H/N_1,\ldots, H/N_t\}\subseteq\tau\mbox{form\,}\mathfrak{X}.
$$

Now let $n>0$. Suppose that $R$ is not a $\sigma$-primary group. Then $A\in\mathfrak{M}$ by Lemma~\ref{l23}. A contradiction. Hence, $R$ is a $\sigma_j$-group for some $\sigma_j\in\sigma\setminus\{\sigma_i\}$. Then $\mbox{O}_{\sigma_j'}(A)=1$ and $F_{\{\sigma_j\}}(A)=\mbox{O}_{\sigma_j}(A)$.
 
Let $f$ be the smallest $l^{\tau}_{\sigma_{n-1}}$-valued $\sigma$-local definition of $\mathfrak{F}=l^{\tau}_{\sigma_n}\mbox{form\,}\mathfrak{M}$ and $h$ be the smallest $l^{\tau}_{\sigma_{n-1}}$-valued $\sigma$-local definition of $\mathfrak{H}=l^{\tau}_{\sigma_n}\mbox{form\,}\mathfrak{X}$. Since for any group $G$ we have
$$
G/F_{\{\sigma_j\}}(G)\simeq (G/\mbox{O}_{\sigma_i}(G))/(F_{\{\sigma_j\}}(G)/\mbox{O}_{\sigma_i}(G))=(G/\mbox{O}_{\sigma_i}(G))/F_{\{\sigma_j\}}(G/\mbox{O}_{\sigma_i}(G)),
$$
in view of Theorem~\ref{t1} we have
$$
f(\sigma_j)=l^{\tau}_{\sigma_{n-1}}\mbox{form\,}(G/F_{\{\sigma_j\}}(G)\mid G\in \mathfrak M)=
$$
$$
=l^{\tau}_{\sigma_{n-1}}\mbox{form\,}((G/\mbox{O}_{\sigma_i}(G))/F_{\{\sigma_j\}}(G/\mbox{O}_{\sigma_i}(G)) \mid G/\mbox{O}_{\sigma_i}(G)\in\mathfrak X)= h(\sigma_j).
$$ 

Since $A\in\mathfrak{F}$, $A/F_{\{\sigma_j\}}(A)\in f(\sigma_j)=h(\sigma_j)\subseteq\mathfrak H$. 
Therefore, $A/\mbox{O}_{\sigma_j}(A)=A/F_{\{\sigma_j\}}(A)\in \mathfrak H$. 
Let $\sigma_a\in\sigma(A)\setminus\{\sigma_j\}$. Then $\mbox{O}_{\sigma_j}(A)\subseteq F_{\{\sigma_a\}}(A)$ and 
$$
A/F_{\{\sigma_a\}}(A)\simeq (A/\mbox{O}_{\sigma_j}(G))/(F_{\{\sigma_a\}}(A)/\mbox{O}_{\sigma_j}(A))=$$
$$(A/\mbox{O}_{\sigma_j}(G))/F_{\{\sigma_a\}}(A/\mbox{O}_{\sigma_j}(A)\in h(\sigma_a).
$$
Thus, for all $\sigma_a\in\sigma(A)$ we have $A/F_{\{\sigma_a\}}(A)\in h(\sigma_a)$.  
Therefore, $A\in\mathfrak{H}$.

Consider now the case when $\mbox{Soc}(A)=N_1\times\ldots\times N_t$, where $t>1$.
Let $M_k$ be the largest normal subgroup in $A$ containing $N_1\times\ldots\times N_{k-1}\times N_{k+1}\times\ldots\times N_t$, but not containing $N_k, k=1,\ldots, t$. Then, in view of Lemma~\ref{l24}, we have $A\in \mbox{R}_0(A/M_1,\ldots, A/M_t)$, where $A/M_k$ is a monolithic group with $\mbox{O}_{\sigma_i}(A/M_k)=1$. It is clear that $A/M_k\in l^{\tau}_{\sigma_n}\mbox{form\,}\mathfrak{M}$. Hence, according to what was proved above, the monolithic group $A/M_k$ belongs to $\mathfrak{H}$ for all $k=1,\ldots, t$. Therefore, $A\in\mathfrak{H}$.
The lemma is proved.

Let $\theta$ be a complete lattice of formations. Then 
following A.N.~Skiba \cite[p. 151]{sk1}, a complete lattice $\theta^\sigma$ of formations we call {\sl $\sigma$-inductive}, if for any set $\{\mathfrak{F}_j \mid j\in J\}\subseteq \theta^\sigma$ and for any set $\{f_j\mid j\in J\}$, where $f_j$ is the integrated  $\theta$-valued $\sigma$-local definition of $\mathfrak{F}_j$, we have
$$
\vee_{\theta^\sigma}(\mathfrak{F}_j\mid j\in j)=LF_\sigma(\vee_{\theta}(f_ji\mid j\in J)).
$$

In particular, if $\sigma=\sigma^1$ we obtain the notion of an inductive lattice of formations \cite[p. 151]{sk1}. 

\begin{theorem} \label{t6} 
The lattice $l^{\tau}_{\sigma_n}$ is $\sigma$-inductive.
\end{theorem}

Proof.
Let $\{\mathfrak{F}_j\mid j\in J\}$ be some collection of $\tau$-closed $n$-multiply $\sigma$-local formations and $f_j$ be some integrate $l^{\tau}_{\sigma_{n-1}}$-valued $\sigma$-local definition of  $\mathfrak{F}_j$.
Let $\mathfrak{F}=\vee_{\sigma_n}^{\tau}(\mathfrak{F}_j\mid j\in J)$, $\mathfrak{M}=LF_\sigma (\vee_{\sigma_{n-1}}^{\tau}(f_j\mid j\in J))$, and $h_j$ be the smallest $l^{\tau}_{\sigma_{n-1}}$-valued $\sigma$-local definition of $\mathfrak{F}_j$.
Then, in view of Lemma~\ref{l21} we have 
$h=\vee_{\sigma_{n-1}}^{\tau}(h_j\mid j\in J)$ is the smallest $l_{\sigma_{n-1}}^{\tau}$-valued $\sigma$-local definition of $\mathfrak{F}$ and, obviously, $h\le f=\vee_{\sigma_{n-1}}^{\tau}(f_j\mid j\in J)$.
Hence, $\mathfrak{F}\subseteq\mathfrak{M}$. Suppose that the reverse inclusion is not true and let $G$ be a group of minimal order in $\mathfrak{M}\setminus\mathfrak{F}$. Then $G$ is a monolithic group with monolith $N=G^{\mathfrak{F}}$.
 
Let $\sigma_i\in\sigma(N)$. 
Suppose that $N$ is not a $\sigma$-primary group. Then $ F_{\{\sigma_i\}}(G)=1$. Therefore,
$$
G\simeq G/ F_{\{\sigma_i\}}(G)\in(\vee_{\sigma_{n-1}}^{\tau}(f_j\mid j\in J))(\sigma_i)=l^{\tau}_{\sigma_{n-1}}\mbox{form\,}(\cup_{j\in J}f_j(\sigma_i)).
$$

Hence, according to Lemma~\ref{l23} we have $ G\in\cup_{j\in J}f_j(\sigma_i)\subseteq\cup_{j\in J}\mathfrak{F}_j\subseteq\mathfrak{F}.$
Contradiction. Therefore, $N$ is a $\sigma_i$-group.
Then $\mbox{O}_{\sigma_i'}(G)=1$ and $F_{\{\sigma_i\}}(G)=\mbox{O}_{\sigma_i}(G)$. Since $G\in\mathfrak{M}=LF_\sigma (\vee_{\sigma_{n-1}}^{\tau}(f_j\mid j\in J))$, we have 
$$
G/\mbox{O}_{\sigma_i}(G)=G/F_{\{\sigma_i\}}(G)\in \vee_{\sigma_{n-1}}^{\tau}(f_j\mid j\in J))(\sigma_i)=l^{\tau}_{\sigma_{n-1}}\mbox{form\,}(\cup_{j\in J}f_j(\sigma_i)).
$$

Hence, by Lemma~\ref{l25} and Theorem~\ref{t1}
$$
G/\mbox{O}_{\sigma_i}(G)\in l^{\tau}_{\sigma_{n-1}}\mbox{form\,}(A/\mbox{O}_{\sigma_i}(A)\mid A\in\cup_{j\in J}f_j(\sigma_i))=
$$
$$
=l^{\tau}_{\sigma_{n-1}}\mbox{form\,}(\cup_{j\in J}l_{\sigma_{n-1}}^{\tau}\mbox{form\,}(A/\mbox{O}_{\sigma_i}(A)\mid A\in f_j(\sigma_i)))=
$$
$$
=l_{\sigma_{n-1}}^{\tau}\mbox{form\,}(\cup_{j\in J}h_j(\sigma_i))=(\vee_{\sigma_{n-1}}^{\tau}(h_j\mid j\in J))(\sigma_i)=h(\sigma_i).
$$
But then $G\in\mathfrak{F}$ in view of Lemma~\ref{l1}. The resulting contradiction shows that $\mathfrak{M}\subseteq\mathfrak{F}$. Therefore, $\mathfrak{F}=\mathfrak{M}$.
The theorem is proved.

In the case when $\tau$ is a trivial subgroup functor from Theorem~\ref{t6} we obtain

\begin{corollary} {\rm \cite{voridr}} \label{c22t6} 
The lattice $l^{\sigma}_n$ is $\sigma$-inductive.
\end{corollary}

In the classical case, when $\sigma=\sigma^1$ from Theorem~\ref{t6} we obtain the following well-known result

\begin{corollary} {\rm \cite[c. 152]{sk1}} \label{c23t6} 
The lattice $l^{\tau}_{n}$ is inductive.
\end{corollary}

It is easy to see that if $\theta$ is a $\sigma$-inductive lattice of formations, then each of its complete sublattices is also $\sigma$-inductive. Thus, in view of Theorem~\ref{t6} we have 

\begin{corollary} \label{c24t6} 
Every complete sublattice of the lattice $l^{\tau}_{\sigma_n}$  is $\sigma$-inductive.
\end{corollary}

\begin{lemma} {\rm \cite[p. 177]{dh}} \label{l26} 
If abelian minimal normal subgroups of $G$ are pairwise non-isomorphic as $G$-modules, then $G$ has a faithful irreducible representation over any field of characteristic 0 or characteristic not dividing $| F(G) |$. 
\end{lemma}

\begin{lemma} \label{l27} 
Let $\mathfrak F$ be a nonempty formation, $\sigma(\mathfrak F)\subseteq\Pi\subseteq\sigma$, and let $f$ be a formation $\sigma$-function such that $f(\sigma_i )=\mathfrak F$ for all $\sigma_i\in\Pi$ and $f(\sigma_i )=\emptyset$ for all $\sigma_i\in\Pi'$. Then $LF_{\sigma } (f)=\mathfrak N_{\Pi} \mathfrak F$ and $f$ is an integrated $\sigma$-local definition of \ $\mathfrak N_{\Pi} \mathfrak F.$ Moreover, if \ $\mathfrak F\in l^{\tau}_{\sigma_{n-1}}$, then $LF_{\sigma } (f)\in l^{\tau}_{\sigma_{n}}$ and $f$ is an integrated $l^{\tau}_{\sigma_{n-1}}$-valued $\sigma$-local definition of\ $\mathfrak N_{\Pi} \mathfrak F.$
\end{lemma}

Proof. Note that if $G\in\mathfrak F$, then $G/F_{\{\sigma_i\}}(G)\in \mathfrak F=f(\sigma_i)$ for all $\sigma_i\in\sigma(G)$ since $\sigma(\mathfrak F)\subseteq\Pi.$ Therefore, $G\in LF_{\sigma } (f)$ and $\mathfrak F\subseteq LF_{\sigma } (f)$, so $f$ is a integrated formation $\sigma$-function. 

 First assume that $LF_{\sigma}(f)\not\subseteq \mathfrak N_{\Pi} \mathfrak F$ and let $G$ be a group of minimal order in $ LF_{\sigma}(f)\setminus\mathfrak N_{\Pi} \mathfrak F$.
Then $R=G^{\mathfrak N_{\Pi} \mathfrak F}$ is an unique minimal normal subgroup of $G$ since $\mathfrak N_{\Pi} \mathfrak F$ is a formation. Let $\sigma_i\in\sigma(R)$.  Since $G\in LF_{\sigma}(f)$, we have $G/F_{\{\sigma_i\}}(G)\in f(\sigma_i)=\mathfrak F$.
If $R$ is not a $\sigma$-primary group, then $F_{\{\sigma_i\}}(G)=1.$ Therefore, $G\simeq G/F_{\{\sigma_i\}}(G)\in f(\sigma_i)=\mathfrak F\subseteq \mathfrak N_{\Pi} \mathfrak F$, a contadiction.  Hence, $R$ is $\sigma$-primary, so $R$ is a $\sigma_i$-group. But then $F_{\{\sigma_i\}}(G)=O_{\sigma_i}(G)$. Therefore, $G/O_{\sigma_i}(G)=G/F_{\{\sigma_i\}}(G)\in f(\sigma_i).$ Thus, $G\in  \mathfrak G_{\sigma_i} f(\sigma_i) =  \mathfrak G_{\sigma_i} \mathfrak F\subseteq \mathfrak N_{\Pi} \mathfrak F.$ This contradiction shows that $LF_{\sigma}(f)\subseteq \mathfrak N_{\Pi} \mathfrak F$.

Now we show that $\mathfrak N_{\Pi} \mathfrak F\subseteq LF_{\sigma}(f).$ Suppose that this is false and let $G$ be a group of minimal order in $\mathfrak N_{\Pi} \mathfrak F\setminus LF_{\sigma}(f)$. Then $G$ has an unique minimal normal subgroup  $R=G^{LF_{\sigma}(f)}.$ In view of the above proved $\mathfrak F\subseteq LF_{\sigma } (f)$, so $G^{\mathfrak F}\ne 1$. Since $G\in \mathfrak N_{\Pi} \mathfrak F,$ $R$ is a $\sigma_i$-group for some $\sigma_i\in\Pi.$ Hence, $G\in \mathfrak G_{\sigma_i} \mathfrak F = \mathfrak G_{\sigma_i}f(\sigma_i )\subseteq LF_{\sigma}(f)$, since $f$ is integrated. A contradiction. Thus, 
$LF_{\sigma}(f)=\mathfrak N_{\Pi} \mathfrak F$.

Finally, let $\mathfrak F\in l^{\tau}_{\sigma_{n-1}}$. Then $f(\sigma_i)\in l^{\tau}_{\sigma_{n-1}}$ for all $\sigma_i\in\mbox{Supp}(f)$. Hence, $f$ is an integrated $l^{\tau}_{\sigma_{n-1}}$-valued $\sigma$-local definition of $\mathfrak N_{\Pi} \mathfrak F.$ But then $LF_{\sigma}(f)\in l^{\tau}_{\sigma_{n}}$ by definition.
The lemma is proved.

\begin{lemma}\label{l28} 
Let $\mathfrak{M}$ and $\mathfrak{H}$ be some $\tau$-closed $n$-multiply $\sigma$-local formations, $\sigma(\mathfrak{M}\cup\mathfrak{H})\subseteq \Pi\subseteq\sigma$. Then  
$$
(\mathfrak{N}_{\Pi}\mathfrak{M})\vee^{\tau}_{\sigma_{n+1}}(\mathfrak{N}_{\Pi}\mathfrak{H})=\mathfrak{N}_{\Pi}(\mathfrak{M}\vee^{\tau}_{\sigma_n}\mathfrak{H}).
$$
In particular, 
$$
(\mathfrak{N}_{\sigma}\mathfrak{M})\vee^{\tau}_{\sigma_{n+1}}(\mathfrak{N}_{\sigma}\mathfrak{H})=
\mathfrak{N}_{\sigma}(\mathfrak{M}\vee^{\tau}_{\sigma_n}\mathfrak{H}).
$$
\end{lemma}

Proof. Let $\mathfrak{X}=\mathfrak{N}_{\Pi}\mathfrak{M}$, $\mathfrak{L}=\mathfrak{N}_{\Pi}\mathfrak{H}$ and $\mathfrak{F}=\mathfrak{N}_{\Pi}(\mathfrak{M}\vee^{\tau}_{\sigma_n}\mathfrak{H})$.
In view of Lemma~\ref{l27} the formations $\mathfrak{X}$, $\mathfrak{L}$ and $\mathfrak{F}$ have the interated $l_{\sigma_n}^{\tau}$-valued $\sigma$-local definitions of $x$, $l$ and $f$ respectively, such that $x(\sigma_i)=\mathfrak{M}$, $l(\sigma_i)=\mathfrak{H}$ and $f(\sigma_i)=\mathfrak{M}\vee^{\tau}_{\sigma_n}\mathfrak{H}$ for all $\sigma_i\in\Pi$ and $x(\sigma_i)=l(\sigma_i)=f(\sigma_i)=\emptyset$ for all $\sigma_i\in\Pi'$. On the other hand, by Theorem~\ref{t6}  we have 
$$
(\mathfrak{N}_{\Pi}\mathfrak{M})\vee^{\tau}_{\sigma_{n+1}}(\mathfrak{N}_{\Pi}\mathfrak{H})=LF_\sigma (x\vee^{\tau}_{\sigma_n}l). 
$$
Therefore, for all $\sigma_i\in\Pi$ we have
$$
(x\vee^{\tau}_{\sigma_n}l)(\sigma_i)=l^{\tau}_{\sigma_n}\mbox{form\,}(x(\sigma_i)\cup l(\sigma_i))=l^{\tau}_{\sigma_n}\mbox{form\,}(\mathfrak{M}\cup\mathfrak{H})=\mathfrak{M}\vee^{\tau}_{\sigma_n}\mathfrak{H}=f(\sigma_i)
$$
and $(x\vee^{\tau}_{\sigma_n}l)(\sigma_i)=\emptyset =f(\sigma_i)$ for all $\sigma_i\in\Pi'$.
Hence, $m\vee^{\tau}_{\sigma_n}h=f$. The latter entails
$$
(\mathfrak{N}_{\Pi}\mathfrak{M})\vee^{\tau}_{\sigma_{n+1}}(\mathfrak{N}_{\Pi}\mathfrak{H})=\mathfrak{N}_{\Pi}(\mathfrak{M}\vee^{\tau}_{\sigma_n}\mathfrak{H}).
$$
In particular, if $\Pi=\sigma$, then $\mathfrak{N}_{\Pi}=\mathfrak{N}_{\sigma}$ and we have
$$
(\mathfrak{N}_{\sigma}\mathfrak{M})\vee^{\tau}_{\sigma_{n+1}}(\mathfrak{N}_{\sigma}\mathfrak{H})=\mathfrak{N}_{\sigma}(\mathfrak{M}\vee^{\tau}_{\sigma_n}\mathfrak{H}).
$$
The lemma is proved.

\begin{lemma} \label{l29}  
Let $\mathfrak{H}$ be some nonempty formation, $\sigma(\mathfrak{H})\subseteq \Pi\subseteq\sigma$, where $|\Pi| > 1$.
Then if \  $\mathfrak{N}_{\Pi}\mathfrak{H}$ is an $n$-multiply $\sigma$-local formation $(n\ge 1)$, then 
$\mathfrak{H}$ is $(n-1)$-multiply $\sigma$-local. In particular, if \ $\mathfrak{N}_{\sigma}\mathfrak{H}\in l^{\sigma}_n$ $(n\ge 1)$, then $\mathfrak{H}\in l^{\sigma}_n$.
\end{lemma} 

Proof. We use induction on $n$. Let $n=1$. Then, since every formation is $0$-multiply $\sigma$-local by definition, the lemma is true.

Let now $n>1$ and the lemma is true for any natural number less than $n$. Then, by Lemma~\ref{l27} the formation $\mathfrak{F}=\mathfrak{N}_{\Pi}\mathfrak{H}$ has an integrated $\sigma$-local definition $f$, such that $f(\sigma_i)=\mathfrak{H}$ for any $\sigma_i\in\Pi$  and $f(\sigma_i)=\emptyset$ for all $\sigma_i\in\Pi'$. 
Let $F$ is the canonical $\sigma$-local definition of $\mathfrak{F}$. Then, in view of Corollary~\ref{c11t2} we have $F(\sigma_i)\in l^{\sigma}_{n-1}$. By Lemma~\ref{l6}~(5) we have $F(\sigma_i)=\mathfrak{G}_{\sigma_i} f(\sigma_i)=\mathfrak{G}_{\sigma_i}\mathfrak{H}$ for all $\sigma_i\in\Pi$ and $F(\sigma_i)=\emptyset$ for all $\sigma_i\in\Pi'$. 
Therefore, for all $\sigma_i\in\Pi$ the formation $\mathfrak{G}_{\sigma_i}\mathfrak{H}$ is $(n-1)$-multiply $\sigma$-local.

Now let $\sigma_i,\sigma_j\in\Pi$, $i\ne j$. Since 
the formations $\mathfrak{G}_{\sigma_i}\mathfrak{H}$ and $\mathfrak{G}_{\sigma_j}\mathfrak{H}$ are $(n-1)$-multiply $\sigma$-local, we have $$
\mathfrak{G}_{\sigma_i}\mathfrak{H}\cap\mathfrak{G}_{\sigma_j}\mathfrak{H}= (\mathfrak{G}_{\sigma_i}\cap\mathfrak{G}_{\sigma_j})\mathfrak{H} = \mathfrak H
$$
is $(n-1)$-multiply $\sigma$-local by Lemma~\ref{l19}. 
In particular, if $\Pi=\sigma$, then $\mathfrak N_\Pi=\mathfrak N_\sigma$  and we obtain the second statement of the lemma. The lemma is proved.

\begin{theorem} \label{t7} 
Let $m$ and $n$ be non-negative integers, where $m>n$. Then the lattice $l^{\tau}_{\sigma_m}$ is not a sublattice of the lattice $l^{\tau}_{\sigma_n}$.
\end{theorem}

Proof. First note that for any formations $\mathfrak{M},$ $\mathfrak{H}\in l^{\tau}_{\sigma_n}$ we have
$$
\mathfrak{M}\vee^{\tau}_{\sigma_n}\mathfrak{H}\subseteq\mathfrak{M}\vee_{\sigma_{n+1}}^{\tau}\mathfrak{H}\subseteq\ldots\subseteq\mathfrak{M}\vee^{\tau}_{\sigma_{m}}\mathfrak{H}.
$$
Therefore, to prove the assertion of the theorem, it suffices to restrict ourselves to the case when $m=n+1$.

We prove this theorem by induction on $n$. Let $n=0$. Denote $\mathfrak{F}=\mathfrak{G}_{\sigma_i}\mathfrak{G}_{\sigma_k}\vee\mathfrak{G}_{\sigma_i}\mathfrak{G}_{\sigma_j}$, where $i$, $k$ and $j$ are pairwise distinct numbers. Let us show that the formation $\mathfrak{F}$ is not $\sigma$-local. Suppose the opposite and let $f$ be the smallest $\sigma$-local definition of $\mathfrak{F}$. 
Therefore, by Lemma~\ref{l21} we have $f=f_1\vee f_2$, where $f_1$ and $f_2$ are the smallest $\sigma$-local definitions of formations $\mathfrak{F}_1=\mathfrak{G}_{\sigma_i}\mathfrak{G}_{\sigma_k}$ and $\mathfrak{F}_2=\mathfrak{G}_{\sigma_i}\mathfrak{G}_{\sigma_j}$ respectively. Then, in view of Theorem~\ref{t1}, we have $f_1(\sigma_i)=\mathfrak{G}_{\sigma_k}$ and $f_2(\sigma_i)=\mathfrak{G}_{\sigma_j}$. Therefore,
$$
f(\sigma_i)=(f_1\vee f_2)(\sigma_i)= f_1(\sigma_i)\vee f_2(\sigma_i)=\mathfrak{G}_{\sigma_k}\vee \mathfrak{G}_{\sigma_j}=\mathfrak{N}_{\{\sigma_k, \sigma_j\}}.
$$
Let $R$ and $Q$ are some groups of order $r\in\sigma_k$ and $q\in\sigma_j$ respectively, $V=R\times Q$ and $p\in\sigma_i$. 
By Lemma~\ref{l26} the group $V$ has a faithful irreducible module $P$ over the field $F_p$ of characteristic $p$.
Let $G=P\rtimes V$. Then since, obviously, $\mbox{O}_{\sigma_i}(G)=P$ and $G/P\simeq R\times Q\in f(\sigma_i)\subseteq\mathfrak F$,  $G\in\mathfrak{F}$ in view of Lemma~\ref{l1}. Note also that $G\notin\mathfrak{G}_{\sigma_i}\mathfrak{G}_{\sigma_k}\cup\mathfrak{G}_{\sigma_i}\mathfrak{G}_{\sigma_j}$. Hence, by Lemma~\ref{l13} 
in $\mathfrak{F}$ there exists a group $H$ with normal subgroups $N,$ $M,$ $N_1,\ldots, N_t,$ $M_1,\ldots, M_t $ $(t\ge 2)$, such that the following statements are hold:
1) $H/N\simeq G$, $M/N=\mbox{Soc}(H/N)$;
2) $N_1\cap\ldots\cap N_t=1$;
3) $H/N_i$ is a monolithic $(\mathfrak{G}_{\sigma_i}\mathfrak{G}_{\sigma_k}\cup\mathfrak{G}_{\sigma_i}\mathfrak{G}_{\sigma_j})$-group with a monolith $M_i/N_i$, which is $H$-somorphic to $M/N$.
Let, for definiteness, $H/N_1\in\mathfrak{G}_{\sigma_i}\mathfrak{G}_{\sigma_k}$. 
Since $\mbox{C}_G(P)=P$, we have $M=\mbox{C}_H(M/N)$. Therefore, $M_1\subseteq M$.
Hence, $V=R\times Q\in\mathfrak{G}_{\sigma_i}\mathfrak{G}_{\sigma_k}$.
A contradiction. Therefore, $\mathfrak{F}$ is not $\sigma$-local. Since the formations $\mathfrak{G}_{\sigma_i},$ $\mathfrak{G}_{\sigma_k}$ and $\mathfrak{G}_{\sigma_j}$ are hereditary and $\sigma$-local,  the formations $\mathfrak{G}_{\sigma_i}\mathfrak{G}_{\sigma_k}$ and $\mathfrak{G}_{\sigma_i}\mathfrak{G}_{\sigma_j}$ are $\tau$-closed by Lemma~\ref{l8} and are $\sigma$-local by Lemma~\ref{l9}. Hence, $\mathfrak{G}_{\sigma_i}\mathfrak{G}_{\sigma_k}$ and $\mathfrak{G}_{\sigma_i}\mathfrak{G}_{\sigma_j}$ belong to the lattice $l^\tau_{\sigma_1}$. However, by Lemma~\ref{l16} we have
$$
\mathfrak{G}_{\sigma_i}\mathfrak{G}_{\sigma_k}\vee^{\tau}\mathfrak{G}_{\sigma_i}\mathfrak{G}_{\sigma_j}=\tau\mbox{form\,}(\mathfrak{G}_{\sigma_i}\mathfrak{G}_{\sigma_k}\cup\mathfrak{G}_{\sigma_i}\mathfrak{G}_{\sigma_j})=
$$
$$
\mbox{form\,}(\mathfrak{G}_{\sigma_i}\mathfrak{G}_{\sigma_k}\cup\mathfrak{G}_{\sigma_i}\mathfrak{G}_{\sigma_j})=
\mathfrak{G}_{\sigma_i}\mathfrak{G}_{\sigma_k}\vee\mathfrak{G}_{\sigma_i}\mathfrak{G}_{\sigma_j},
$$
and, as shown above, $\mathfrak{G}_{\sigma_i}\mathfrak{G}_{\sigma_k}\vee\mathfrak{G}_{\sigma_i}\mathfrak{G}_{\sigma_j}\notin l^\tau_{\sigma_1}$. Therefore, the lattice $l^{\tau}_{\sigma_1}$ is not a sublattice of the lattice $l^{\tau}_{\sigma_0}$.

Now let $n>1$ and suppose that the theorem is true for $n-1$. Then there are $\tau$-closed $n$-multiply $\sigma$-local formations $\mathfrak{M}$ and $\mathfrak{H}$, such that $\mathfrak{M}\vee^{\tau}_{\sigma_{n-1}}\mathfrak{H}\notin l^{\tau}_{\sigma_n}$. Let $\mathfrak{M}_1=\mathfrak{N}_{\sigma}\mathfrak{M}$, $\mathfrak{H}_1=\mathfrak{N}_{\sigma}\mathfrak{H}$.
By Lemma~\ref{l27} the formations $\mathfrak{M}_1$ and $\mathfrak{H}_1$ have an interated $l^{\tau}_{\sigma_n}$-valued $\sigma$-local definitions $m$ and $h$ and $l^{\tau}_{\sigma_{n-1}}$-valued $\sigma$-local definitions $m_1$ and $h\,_1$ respectively, such that $m(\sigma_i)=m_1(\sigma_i)=\mathfrak{M}$, $h(\sigma_i)=h\,_1(\sigma_i)=\mathfrak{H}$ for all $\sigma_i\in\sigma$. Therefore, each of the formations $\mathfrak{M}_1$ and $\mathfrak{H}_1$ belongs to the lattice $l^{\tau}_{\sigma_{n+1}}$.
Suppose $\mathfrak{M}_1\vee^{\tau}_{\sigma_n}\mathfrak{H}_1\in l^{\tau}_{\sigma_{n+1}}$.
Then, by Lemma~\ref{l28} we have
$
\mathfrak{M}_1\vee^{\tau}_{\sigma_n}\mathfrak{H}_1=\mathfrak{N}_{\sigma}(\mathfrak{M}\vee^{\tau}_{\sigma_{n-1}}\mathfrak{H})
$
and
$$
\mathfrak{M}_1\vee^{\tau}_{\sigma_n}\mathfrak{H}_1=l_{\sigma_n}^{\tau}\mbox{form\,}(\mathfrak{M}_1\cup\mathfrak{H}_1)=l^{\tau}_{\sigma_{n+1}}\mbox{form\,}(\mathfrak{M}_1\cup\mathfrak{H}_1)=\mathfrak{M}_1\vee_{\sigma_{n+1}}^{\tau}\mathfrak{H}_1.
$$
Hence, $\mathfrak{M}\vee^{\tau}_{\sigma_{n-1}}\mathfrak{H}$ is $n$-multiply $\sigma$-local  by Lemma~\ref{l29}. Therefore, $\mathfrak{M}\vee^{\tau}_{\sigma_{n-1}}\mathfrak{H}\in l^{\tau}_{\sigma_n}$. The latter contradicts the choice of the formations $\mathfrak{M}$ and $\mathfrak{H}$. Thus, the lattice $l^{\tau}_{\sigma_{n+1}}$ is not a sublattice of the lattice $l^{\tau}_{\sigma_n}$. The theorem is proved.

In the case when $\sigma =\sigma ^{1}$, we obtain the following well-known result from Theorem~\ref{t7}.

\begin{corollary} {\rm \cite[p. 157]{sk1}} \label{c25t7} 
For any non-negative integers $m$ and $n$, where $m>\, n$, the lattice $l^{\tau}_{m}$ is not a sublattice in $l^{\tau}_{n}$.
 \end{corollary}
 
At the same time, the following theorem holds.

\begin{theorem} \label{t8} 
The lattice $l^{\tau}_{\sigma_n}$ is a complete sublattice of the lattice $l^\sigma_n$.
\end{theorem}

Proof. We use induction on $n$. For $n=0$ the theorem is true in view of Lemma~\ref{l16}.
Let $n>0$ and assume that for $n-1$ the theorem is true. Let $\{\mathfrak{F}_j\mid j\in J\}$ be an arbitrary set of $\tau$-closed $n$-multiply $\sigma$-local formations $\mathfrak{F}_j=LF_{\sigma}(f_j)$, where $f_j$ is the smallest $l^{\tau}_{n-1}$-valued $\sigma$-local  definition of $\mathfrak{F}_j$. By Theorem~\ref{t6} for the trivial subgroup functor $\tau$ we have
$$
\vee^\sigma_n(\mathfrak{F}_j\mid j\in J)=LF(\vee^\sigma_{n-1}(f_j\mid j\in J)).
$$
By the inductive assumption, for any $\sigma_i\in\sigma(\cup_{j\in J}\mathfrak{F}_j)$ the formation
$$
(\vee^\sigma_{n-1}(f_j\mid j\in J))(\sigma_i)=\vee_{\sigma_{n-1}}^{\tau}(\mathfrak{F}_j(\sigma_i)\mid j\in J)
$$
is $\tau$-closed. Therefore, the formation $\vee^\sigma_n(\mathfrak {F}_j\mid j\in J)=LF_{\sigma}(h)$, where $h=\vee^\sigma_{n-1}(f_j\mid j\in J)$, is $\tau$-closed by Theorem~\ref{t3}. Thus, the lattice $l^{\tau}_{\sigma_n}$ is a complete sublattice of the lattice $l^\sigma_n$. The theorem is proved.

\begin{corollary} \label{c26t8}  
The lattice $l^{\tau}_{\sigma}$ is a complete sublattice of the lattice $l_\sigma$.
\end{corollary} 
 
In the case $\sigma =\sigma ^{1}$ from Theorem~\ref{t8} we obtain the following well-known result.

\begin{corollary} {\rm \cite[p. 158]{sk1}} \label{c27t8} 
The lattice $l^{\tau}_{n}$ is a complete sublattice of the lattice $l_n$.
\end{corollary}

\begin{lemma} \label{l30} 
Let $\mathfrak F=l^{\tau}_{\sigma_n}{\rm form\,}(\mathfrak X)=LF_{\sigma } (f)$, where $f$ is a  $l^{\tau}_{\sigma_{n-1}}$-valued $\sigma$-local definition of $\mathfrak F.$ Let $h$ be a formation $\sigma$-function such that $h(\sigma_i )=l^{\tau}_{\sigma_{n-1}}{\rm form\,}(\mathfrak X(\sigma_i ))$ for all $\sigma_i \in \Pi =\sigma (\mathfrak X)$ and $h(\sigma_i )=\emptyset $ for all $\sigma_i \in \Pi '.$
Then $\Pi =\sigma (\mathfrak F),$ $h$ is $l^{\tau}_{\sigma_{n-1}}$-valued $\sigma$-local definition of $\mathfrak F$ and $h(\sigma_i )\subseteq f(\sigma_i )\cap \mathfrak F$ for all $i.$
\end{lemma}

Proof. In view of Theorem~\ref{t1} we have $\Pi =\sigma (\mathfrak F).$  
Let $\mathfrak H=LF_{\sigma } (h).$ Then, obviously, $\mathfrak X\subseteq \mathfrak H.$ On the other hand, since $h$ is a $l^{\tau}_{\sigma_{n-1}}$-valued formation $\sigma$-function, then $\mathfrak H$ is $l^{\tau}_{\sigma_{n}}$-formation.
Therefore, $\mathfrak F\subseteq \mathfrak H.$
Since $\mathfrak X(\sigma_i )\subseteq f(\sigma_i )$ and $f(\sigma_i )$ is a $l^{\tau}_{\sigma_{n-1}}$-formation, then $h(\sigma_i )\subseteq f(\sigma_i )\cap \mathfrak F$ for all $\sigma_i \in \sigma.$ Therefore, $\mathfrak H\subseteq \mathfrak F$ and $\mathfrak F=\mathfrak H$. The lemma is proved.

\begin{lemma} \label{l31} 
Let $\mathfrak{F}_j=LF_{\sigma}(f _j )$, where $f _j $ is a $l^{\tau}_{\sigma_{n-1}}$-valued $\sigma$-local definition of $\mathfrak{F}_j$, $j = 1, 2$. 
Then $\mathfrak{F} = \mathfrak{F}_1 \vee^{\tau}_{\sigma_{n}}\mathfrak{F}_2 =LF_{\sigma} (f)$, where $f = f _1 \vee^{\tau}_{\sigma_{n-1}} f _2$ is integrated $l^{\tau}_{\sigma_{n-1}}$-valued $\sigma$-local definition of $\mathfrak{F}$.
\end{lemma}

Proof. Let $h_j$ be the smallest $l^{\sigma}_{n-1}$-valued $\sigma$-local definition of $\mathfrak{F}_j $ and let $F_{j}$ be the canonical $\sigma$-local definition of $\mathfrak{F}_{j}$, $j = 1, 2$. Let $h$ be the smallest $l^{\sigma}_{n-1}$-valued $\sigma$-local definition of $\mathfrak{F}$ and $F$ be the canonical $\sigma$-local definition of $\mathfrak{F}$. Then for each $\sigma_i\in \sigma$ we have $h_j(\sigma_i) \subseteq f _j(\sigma_i) \subseteq F_j (\sigma_i)$ by Lemmas~\ref{l6}~(5) and~\ref{l30}. Moreover, it also follows from Lemmas~\ref{l6}~(5) and~\ref{l30} that   
$$
h(\sigma_i) = l^{\sigma}_{n-1}\text{form\,}((\mathfrak{F}_1 \cup \mathfrak{F}_2)(\sigma_i))=
$$
$$= l^{\sigma}_{n-1}\text{form\,}(\mathfrak{F}_1(\sigma_i) \cup \mathfrak{F}_2(\sigma_i))=l^{\sigma}_{n-1}\text{form\,}(h_1(\sigma_i) \cup h_2(\sigma_i))\subseteq$$
$$ \subseteq f(\sigma_i) \subseteq \mathfrak{G}_{\sigma_i}l^{\sigma}_{n-1}\text{form\,}(h_1(\sigma_i) \cup h_2(\sigma_i)) \subseteq \mathfrak{G}_{\sigma_i} h(\sigma_i) = F(\sigma_i).
$$
Therefore, $ h(\sigma_i) \subseteq f(\sigma_i) \subseteq F(\sigma_i) $ for all $\sigma_i\in \sigma$, so $\mathfrak{F} = LF_{\sigma} (f)$. The lemma is proved.

\begin{theorem} \label{t9} 
The lattice $l^{\tau}_{\sigma_n}$ is modular.
\end{theorem}

{Proof.} We prove this theorem by induction on $n$. If $n=0$, then $l^{\tau}_{\sigma_0}=l^{\tau}_{0}=\tau$ and the theorem is true in view of Corollary~4.2.8 \cite[p. 167]{sk1}. Now let $n > 0$ and assume that the lattice $l^{\tau}_{\sigma_r}$ is modular for all $r < n$. 

Let $\mathfrak{F}_j = LF_{\sigma}(f_j)$ be a $\tau$-closed $n$-multiply $\sigma$-local formation, $j = 1, 2, 3$, where $\mathfrak{F}_2 \subseteq \mathfrak{F}_1$ and let $f_j$ be the smallest $l^{\tau}_{\sigma_{n-1}}$-valued $\sigma$-local definition of $\mathfrak{F}_j $.
Let us show that 
$$ \mathfrak{F}_1 \cap (\mathfrak{F}_2\vee^{\tau}_{\sigma_n} \mathfrak{F}_3) = \mathfrak{F}_2 \vee^{\tau}_{\sigma_n} (\mathfrak{F}_1\cap \mathfrak{F}_3).
$$ 

Note that the formation $\sigma$-functions $f_1,$ $f_2$ and $ f_3$ are integrated, and $f _2 (\sigma_i) \subseteq f _1 (\sigma_i)$ for all $i $ by Corollary~\ref{c2t1}. Then $f _2 \vee^{\tau}_{\sigma_{n-1}} f _3$ is integrated $l^{\tau}_{\sigma_{n-1}}$-valued and $\mathfrak{F}_2 \vee^{\tau}_{\sigma_n} \mathfrak{F}_3 = LF_{\sigma} (f _2 \vee^{\tau}_{\sigma_{n-1}} f _3)$ by Lemma~\ref{l31}. Therefore,  in view of Lemma~\ref{l4} we have  $\mathfrak{F}_1 \cap (\mathfrak{F}_2\vee^{\tau}_{\sigma_n} \mathfrak{F}_3)= LF_{\sigma}(f_{1}\cap (f _2 \vee^{\tau}_{\sigma_{n-1}} f _3))$. 

On the other hand, $f_{1}\cap f_{3}$ is an interated $l^{\tau}_{\sigma_{n-1}}$-valued $\sigma$-local definition of $\mathfrak{F}_1\cap \mathfrak{F}_3$ by Lemma~\ref{l4}. Now applying Lemma~\ref{l31} we obtain $\mathfrak{F}_2 \vee^{\tau}_{\sigma_n} (\mathfrak{F}_1\cap \mathfrak{F}_3)=LF_{\sigma }(f_{2}\vee^{\tau}_{\sigma_{n-1}} (f_{1}\cap f_{3})).$
Since the functions $f_1, f_2, f_3$ are $l^{\tau}_{\sigma_{n-1}}$-valued, by the inductive hypothesis we obtain that for each $i$ we have
$$ 
f _1 (\sigma_i) \cap (f _2 (\sigma_i) \vee^{\tau}_{\sigma_{n-1}} f _3(\sigma_i)) = f _2(\sigma_i) \vee^{\tau}_{\sigma_{n-1}} (f _1(\sigma_i) \cap f _3(\sigma_i)).
$$ 
Therefore, $ f_{1}\cap (f _2 \vee^{\tau}_{\sigma_{n-1}} f _3))=f_{2}\vee^{\tau}_{\sigma_{n-1}} (f_{1}\cap f_{3}))$. But then we have 
$$ \mathfrak{F}_1 \cap (\mathfrak{F}_2\vee^{\tau}_{\sigma_n} \mathfrak{F}_3) = \mathfrak{F}_2 \vee^{\tau}_{\sigma_n} (\mathfrak{F}_1\cap \mathfrak{F}_3).
$$ 
The theorem is proved. 

In particular, if $\tau$ is a trivial subgroup functor from Theorem~\ref{t9} we obtain the following result.

\begin{corollary} {\rm \cite{chsafsk}} \label{c28t9} 
The lattice $l^{\sigma}_n$ is modular.
\end{corollary}
 
In the classical case, when $\sigma=\sigma^1$ from Theorem~\ref{t9} we obtain thewell-known results
 
\begin{corollary} {\rm \cite[p. 167]{sk1}} \label{c29t9} 
The lattice $l^\tau_n$ is modular.
\end{corollary}
 
\begin{corollary} {\rm \cite[p. 105]{shsk}} \label{c30t9} 
For any non-negative integer $n$ the lattice $l_n$ is modular.
\end{corollary}

For any two $\tau$-closed $n$-multiply $\sigma$-local formations
 $\mathfrak{M}$ and $\mathfrak{H}$, where $\mathfrak{M}\subseteq\mathfrak{H}$, denote by $\mathfrak{H}/^{\tau}_{\sigma_n}\mathfrak{M}$  the lattice $\tau$-closed $n$-multiply $\sigma$-local formations $\mathfrak X$ such that $\mathfrak{M}\subseteq\mathfrak X\subseteq\mathfrak{H}$. 

\begin{corollary} \label{c31t9} 
For any $\tau$-closed $n$-multiply $\sigma$-local formations $\mathfrak{M}$ and $\mathfrak{H}$ the lattice isomorphism is holds
$$
\mathfrak{M}\vee^{\tau}_{\sigma_n}\mathfrak{H}/^{\tau}_{\sigma_n}\mathfrak{M}\simeq\mathfrak{H}/^{\tau}_{\sigma_n}\mathfrak{H}\cap \mathfrak{M}.
$$
\end{corollary}
 
In particular, if $\tau$ is a trivial subgroup functor, we have
  
\begin{corollary} \label{c32t9} 
For any $n$-multiply $\sigma$-local formations $\mathfrak{M}$ and $\mathfrak{H}$ the lattice isomorphism is holds
 $$
\mathfrak{M}\vee^{\sigma}_n\mathfrak{H}/^{\sigma}_n\mathfrak{M}\simeq\mathfrak{H}/^{\sigma}_n\mathfrak{H}\cap \mathfrak{M}.
$$
\end{corollary}
 
In the case $\sigma=\sigma^1$ we obtain
 
\begin{corollary}  {\rm \cite[p. 168]{sk1}} \label{c33t9} 
For any $\tau$-closed $n$-multiply local formations $\mathfrak{M}$ and $\mathfrak{H}$ the lattice isomorphism is holds
$$
\mathfrak{M}\vee^{\tau}_n\mathfrak{H}/^{\tau}_n\mathfrak{M}\simeq\mathfrak{H}/^{\tau}_n\mathfrak{H}\cap \mathfrak{M}.
$$
\end{corollary}

\begin{corollary}  {\rm \cite[p. 168]{sk1}} \label{c34t9} 
For any $n$-multiply local formations $\mathfrak{M}$ and $\mathfrak{H}$  the lattice isomorphism is holds
$$
\mathfrak{M}\vee_n\mathfrak{H}/_n\mathfrak{M}\simeq\mathfrak{H}/_n\mathfrak{H}\cap \mathfrak{M}.
$$
\end{corollary}

\begin{theorem} \label{t10} 
The lattice $l^{\tau}_{\sigma_n}$ is algebraic.
\end{theorem}

{Proof.} First let us show that every one-generated $\tau$-closed $n$-multiply $\sigma$-local formation $ \mathfrak{F}= l^{\tau}_{\sigma_n} \text{\rm form\,}(G)$ is a compact element in $l^{\tau}_{\sigma_n}$. Let $ \mathfrak{F}\subseteq \mathfrak{M}=l^{\tau}_{\sigma_n} \text{\rm form\,}(\cup_{j\in J}\mathfrak{F}_j)$ for some set $\{\mathfrak{F}_j \mid j\in J \}\subseteq  l^{\tau}_{\sigma_n}$. 
Let $f_j$ be the smallest $l^{\tau}_{\sigma_{n-1}}$-valued $\sigma$-local definition of $\mathfrak{F}_j$ for all $j\in J$, $f$ be the smallest $l^{\tau}_{\sigma_{n-1}}$-valued $\sigma$-local definition of $\mathfrak{F}$ and let $m$ be the smallest $l^{\tau}_{\sigma_{n-1}}$-valued $\sigma$-local definition of $\mathfrak{M}$.  

Let $\sigma_i\in \sigma (G)$. Then $f(\sigma_i)=l^{\tau}_{\sigma_{n-1}}\text{form\,}(G/F_{\{\sigma_i\}}(G))$ and $m(\sigma_i)=l^{\tau}_{\sigma_{n-1}}\text{form\,}(\bigcup_{j\in J}f_{j}(\sigma_i))$ by Theorem~\ref{t1}. In addition, it follows from Corollary~\ref{c2t1} that $f(\sigma_i)\subseteq  m(\sigma_i)$ since $\mathfrak{F}\subseteq \mathfrak{M}$. Since $f(\sigma_i)$  is a one-generated $\tau$-closed $(n-1)$-multiply $\sigma$-local formation, by induction there exist $j_1,\ldots,j_t\in J$ such that 
$$G/F_{\{\sigma_i\}}(G)\in  f(\sigma_i)\subseteq f_{j_1}(\sigma_i)\vee^{\tau}_{\sigma_{n-1}}\ldots\vee^{\tau}_{\sigma_{n-1}}f_{j_t}(\sigma_i).
$$
Since in this case $|\sigma (G)|<\infty$, it follows that for some $j_1,\ldots,j_k\in J$ we have $G\in \mathfrak{F}_{j_1}\vee^{\tau}_{\sigma_n}\cdots\vee^{\tau}_{\sigma_n}\mathfrak{F}_{j_k}$. Therefore, $\mathfrak{F}\subseteq\mathfrak{F}_{j_1}\vee^{\tau}_{\sigma_n} \cdots\vee^{\tau}_{\sigma_n}\mathfrak{F}_{j_k}$. 

It is clear that for any $\tau$-closed $n$-multiply $\sigma$-local formation $ \mathfrak{F}$ we have $\mathfrak{F}=l^{\tau}_{\sigma_n} \text{\rm form\,}(\cup_{k\in K}{\mathfrak{F}}_k),$ where $\{\mathfrak{F}_k  \mid  k\in K\}$ is the set of all one-generated $\tau$-closed $n$-multiply $\sigma$-local formations contained in $\mathfrak{F}$. Therefore, the lattice $l^{\tau}_{\sigma_n}$ is algebraic.
The theorem is proved.

If $\tau$ is a trivial subgroup functor, we have

\begin{corollary} {\rm \cite{chsafsk}} \label{c35t10} 
The lattice $l^{\sigma}_n$ of all $n$-multiply $\sigma$-local formations is algebraic.
\end{corollary}

In the case $\sigma =\sigma ^{1}$ from Theorem~\ref{t10} we obtain the following well-known result.

\begin{corollary} {\rm \cite[p. 179]{sk1}}  \label{c36t10} 
The lattice $l^{\tau}_n$ is algebraic.
\end{corollary}

Let $\mathfrak{X}$ be some nonempty class of groups. The complete lattice of formations $\theta$ is called \cite[p. 159]{sk1} $\mathfrak{X}$-{\sl separable}, if for any term $\nu (x_1,\ldots , x_m)$ signatures $\{\cap,\vee_{\theta}\}$, any $\theta$-formations $\mathfrak{F}_1,\ldots, \mathfrak{F}_m$ and any group $A\in\mathfrak{X}\cap\nu (\mathfrak{F}_1,\ldots ,\mathfrak{F}_m)$ there are $\mathfrak{X}$-groups $A_1,\ldots, A_m$, such that $A\in\nu (\theta\mbox{form\,}A_1, \ldots ,\theta\mbox{form\,}A_m)$.

\begin{theorem} \label{t11} 
{\sl The lattice $l^{\tau}_{\sigma_n}$ is $\mathfrak{G}$-separable}.
\end{theorem}

Proof. Let $\nu (x_1,\ldots , x_m)$ be the term signature $\{\cap,\vee^{\tau}_{\sigma_n}\}$, $\mathfrak{F}_1,\ldots ,\mathfrak{F}_m$ are arbitrary formations from $l^{\tau}_{\sigma_n}$ and $A\in\nu (\mathfrak{F}_1,\ldots ,\mathfrak{F}_m)$.
By induction on the number $k$ of occurrences of symbols from $\{\cap,\vee^{\tau}_{\sigma_n}\}$ into the term $\nu$ we will show that there are groups $A_j\in\mathfrak{F}_j \quad (j=1,\ldots,m)$, such that $A\in\nu (\mathfrak{M}_1,\ldots,\mathfrak{M}_m)$, where $\mathfrak{M}_j=l^{\tau}_{\sigma_n}\mbox{form\,}A_j$.

If $k=0$, then the statement of the theorem is obvious. Let $k=1$. Then $A\in\mathfrak F_1\cap\mathfrak F_2$, or $A\in\mathfrak{F}_1\vee^{\tau}_{\sigma_n}\mathfrak{F}_2$. If $A\in\mathfrak F_1\cap\mathfrak F_2$, then $A\in l^{\tau}_{\sigma_n}\mbox{form\,}A\cap l^{\tau}_{\sigma_n}\mbox{form\,}A$.

Let $A\in\mathfrak{F}_1\vee^{\tau}_{\sigma_n}\mathfrak{F}_2$.
Let us prove the assertion of the theorem by induction on $n$.

Let $n=0$. Then $A\in\mathfrak{F}_1\vee^{\tau}_{\sigma_0}\mathfrak{F}_2=l^{\tau}_{\sigma_0}\mbox{form\,}(\mathfrak{F}_1\cup\mathfrak{F}_2)=\tau \mbox{form\,}(\mathfrak{F}_1\cup\mathfrak{F}_2).$
According to Lemma~\ref{l15} the group $A$ is a homomorphic image of some group $H$, where $H\in \mbox{R}_0\mbox{\rm S}_{\overline\tau}(\mathfrak{F}_1\cup\mathfrak{F}_2)$. By Lemma~\ref{l16} we have $\mbox{R}_0\mbox{\rm S}_{\overline\tau}(\mathfrak{F}_1\cup\mathfrak{F}_2)= \mbox{R}_0(\mathfrak{F}_1\cup\mathfrak{F}_2)$. Hence, $H\in \mbox{R}_0(\mathfrak{F}_1\cup\mathfrak{F}_2)$. 

Let us show that $H^{\mathfrak{F}_1}\cap H^{\mathfrak{F}_2}=1$. 
Indeed, if $H\in\mathfrak F_i,$ then $H^{\mathfrak{F}_i}=1$ and $H^{\mathfrak{F}_i}\cap H^{\mathfrak{F}_j}=1$, $i\ne j , \ i,j\in\{1, 2\}.$ Let $H\notin\mathfrak F_i,$ $i=1, 2.$
Since $H\in \mbox{R}_0(\mathfrak{F}_1\cup\mathfrak{F}_2)$, by the definition of the operation ${\rm R}_0$ the group $H$ has normal subgroups $N_1, \ldots, N_k$, such that $\cap_{s=1}^k N_s=1$ and $H/K_s\in \mathfrak{F}_1\cup\mathfrak{F}_2$,  $s=1,\ldots, k.$
Without loss of generality, we can assume that $H/K_s\in \mathfrak{F}_1$ for $s=1, \ldots, t$ and $H/K_s\in \mathfrak{F}_2$ for $s=t+1, \ldots, k$. Then $H^{\mathfrak F_1}\subseteq \cap_{s=1}^t N_s$ and $H^{\mathfrak F_2}\subseteq \cap_{s=t+1}^k N_s$. Therefore, 
$$
H^{\mathfrak{F}_1}\cap H^{\mathfrak{F}_2}\subseteq (\cap_{s=1}^t N_s)\cap (\cap_{s=t+1}^k N_s) = 1 \text{ and } H\in \mbox{R}_0(H/H^{\mathfrak F_1}, H/H^{\mathfrak F_2}). 
$$
Hence,
$$
A\in\mbox{form\,}(H/H^{\mathfrak{F}_1},H/H^{\mathfrak{F}_2})=\mbox{form\,}(H/H^{\mathfrak{F}_1})\vee\mbox{form\,}(H/H^{\mathfrak{F}_2}) \subseteq \mathfrak{F}_1\vee^{\tau}_{\sigma_0}\mathfrak{F}_2.
$$

Now let $n>0, $ $\sigma (A)=\{\alpha_1,\ldots ,\alpha_t\}$ and $A\in\mathfrak{F}_1\vee^{\tau}_{\sigma_n}\mathfrak{F}_2$.
Denote by $f_j$ the smallest $l^{\tau}_{\sigma_{n-1}}$-valued $\sigma$-local definition of $\mathfrak{F}_j, j=1, 2$. 
Then $f_1\vee^{\tau}_{\sigma_{n-1}}f_2$ is the smallest $l^{\tau}_{\sigma_{n-1}}$-valued $\sigma$-local definition of $\mathfrak{F}_1\vee^{\tau}_{\sigma_n}\mathfrak{F}_2$ in view of Lemma~\ref{l21}. Hence, by Theorem~\ref{t1} 
$$
A/ F_{\{\alpha_i\}}(A)\in f_1(\alpha_i)\vee^{\tau}_{\sigma_{n-1}}f_2(\alpha_i).
$$
By induction, there are groups $A_{i_1}\in f_1(\alpha_i)$, $A_{i_2}\in f_2(\alpha_i)$, such that
$$
A / F_{\{\alpha_i\}}(A)\in(l^{\tau}_{\sigma_{n-1}}\mbox{form\,}A_{i_1})\vee^{\tau}_{\sigma_{n-1}}(l^{\tau}_{\sigma_{n-1}}\mbox{form\,}A_{i_2}).
$$
It is clear that 
$
(l^{\tau}_{\sigma_{n-1}}\mbox{form\,}A_{i_1})\vee^{\tau}_{\sigma_{n-1}}(l^{\tau}_{\sigma_{n-1}}\mbox{form\,}A_{i_2})=l^{\tau}_{\sigma_{n-1}}\mbox{form\,}(A_{i_1}, A_{i_2}).
$

We denote by $\mathfrak{L}_1$ the $\tau$-closed semiformation generated by the group $A_{i_1}$ and by $\mathfrak{L}_2$ the $\tau$-closed semiformation, generated by the group $A_{i_2}$. Then, in view of Lemma~\ref{l14} we have $\mathfrak{L}_1=(A_1,\ldots, A_t)$ and $\mathfrak{L}_2=(B_1,\ldots ,B_r)$ for some $A_1,\ldots, A_t\in \mbox{QS}_{\overline{\tau}}\,(A_{i_1})$ and $B_1,\ldots, B_r\in\mbox{QS}_{\overline{\tau}}\,(A_{i_2})$. Note that $\mathfrak{L}_1\cup\mathfrak{L}_2$ is a $\tau$-closed semiformation as the union of $\tau$-closed semiformations. Then we have
$$
A/ F_{\{\alpha_i\}}(A)\in l^{\tau}_{\sigma_{n-1}}\mbox{form\,}(A_{i_1},A_{i_2})=l_{\sigma_{n-1}}^{\tau}\mbox{form\,}(\mathfrak{L}_1\cup\mathfrak{L}_2)=
$$
$$= l_{\sigma_{n-1}}^{\tau}\mbox{form\,}(A_1,\ldots,A_t, B_1,\ldots ,B_r).
$$
Without loss of generality, in view of Lemma~\ref{l25}, we can assume that  $O_{\alpha_i}(A_k)=1$ and $O_{\alpha_i}(B_l)=1$, for all $k=1,\ldots ,t$ and $l=1,\ldots, r$. Let $D_{i_1}=A_1\times\ldots \times A_t$ and $D_{i_2}=B_1\times\ldots \times B_r$. Then $O_{\alpha_i}(D_{i_1})=1$ and $O_{\alpha_i}(D_{i_2})=1$.
Note also that
$$
A/ F_{\{\alpha_i\}}(A)\in l^{\tau}_{\sigma_{n-1}}\mbox{form\,}(D_{i\,_1},D_{i\,_2}) \subseteq l_{\sigma_{n-1}}^{\tau}\mbox{form\,}(A_{i_1},A_{i_2}).
$$

Let $P_i$ be some non-indentity $\alpha_i$-group. Let us denote by $V_{i_1}$  and $V_{i_2}$ the regular wreath products of the group $P_i$ with the groups $D_{i_1}$ and $D_{i_2}$ respectively. Then $V_{i_1}=P_i\wr D_{i_1}=Y_{i_1}\rtimes D_{i_1}$ and $V_{i_2}=P_i\wr D_{i_2}=Y_{i_2}\rtimes D_{i_2}$, where $Y_{i_1}$ and $Y_{i_2}$ are the base groups of the wreath products $V_{i_1}$ and $V_{i_2}$ respectively. Since $O_{\alpha_i'}(V_{i_1})=1$,  we have $F_{\{\alpha_i\}}(V_{i_1})=O_{\alpha_i}(V_{i_1})$ and 
$$
O_{\alpha_i}(V_{i_1})=O_{\alpha_i}(V_{i_1})\cap Y_{i_1}D_{i_1}=Y_{i_1}(O_{\alpha_i}(V_{i_1})\cap D_{i_1})=Y_{i_1}.
$$
Hence, $F_{\{\alpha_i\}}(V_{i_1})=O_{\alpha_i}(V_{i_1})=Y_{i_1}$. Similarly, $F_{\{\alpha_i\}}(V_{i_2})=O_{\alpha_i}(V_{i_2})=Y_{i_2}$.

Then $V_{i_1}/O_{\alpha_i}(V_{i_1})=V_{i_1}/Y_{i_1}\simeq D_{i_1}\in l^{\tau}_{\sigma_{n-1}}\mbox{form\,}A_{i_1}\subseteq f_{1}(\alpha_i)$. By Lemma~\ref{l1} we have $V_{i_1}\in\mathfrak{F}_1$. Similarly, $V_{i_2}\in\mathfrak{F}_2$. Then
$$
A_1=B_{1_1}\times B_{2_1}\times\ldots\times B_{t_1}\in\mathfrak{F}_1, \  \  A_2=B_{1_2}\times B_{2_2}\times\ldots\times B_{t_2}\in\mathfrak{F}_2.
$$

Let us now show that $A/ F_{\{\alpha_i\}}(A)\in f(\alpha_i)$, where $f$ is the smallest $l^{\tau}_{\sigma_{n-1}}$-valued $\sigma$-local definition of $\mathfrak{F}$.

Since $V_{i_1}\in\mathfrak{F}_1\subseteq\mathfrak F$, we have $V_{i_1}\in\mathfrak{F}$. Hence, $D_{i_1}\simeq V_{i_1}/ F_{\{\alpha_i\}}(V_{i_1})\in f(\alpha_i)$ by Theorem~\ref{t1}.
Similarly, it is proved that $D_{i_2}\in f(\alpha_i)$.
Therefore, the inclusion $l^{\tau}_{\sigma_{n-1}}\mbox{form\,}(D_{i_1}, D_{i_2})\subseteq f(\alpha_i)$ holds. But then
 
$$
A/ F_{\{\alpha_i\}}(A)\in l^{\tau}_{\sigma_{n-1}}\mbox{form\,}(D_{i_1}, D_{i_2})\subseteq f(\alpha_i).
$$

Consequently, $A/ F_{\{\alpha_i\}}(A)\in f(\alpha_i)$ for all $\alpha_i\in\sigma(A)$. Therefore,
$$
A\in\mathfrak{F}=(l^{\tau}_{\sigma_n}\mbox{form\,}A_1)\vee^{\tau}_{\sigma_n}(l^{\tau}_{\sigma_n}\mbox{form\,}A_2).
$$
Thus, for $k=1$ the theorem is true.

Now let the term $\nu$ has $k>1$ occurrences of characters from $\{\cap, \vee^{\tau}_{\sigma_n}\}$  and for terms with fewer occurrences of characters from $\{\cap, \vee^{\tau}_{\sigma_n}\}$ the theorem is true.

Let the term $\nu$ have the form
$$
\nu_1(x_{i_1},\ldots ,x_{i_a})\triangle\,\nu_2(x_{j_1},\ldots, x_{j_b}),
$$
where $\triangle\,\in\{\cap,\vee^{\tau}_{\sigma_n}\}$, and $ \{x_{i_1},\ldots ,x_{i_a}\}\cup\{x_{j_1},\ldots, x_{j_b}\}=\{x_1,\ldots, x_m\}.$
Put $\mathfrak{H}_1=\nu_1(\mathfrak{F}_{i_1},\ldots ,\mathfrak{F}_{i_a})$ and $\mathfrak{H}_2=\nu_2(\mathfrak{F}_{j_1},\ldots ,\mathfrak{F}_{j_b})$.
Then, by what was proved above, there are groups $A_1\in\mathfrak{H}_1$, $A_2\in\mathfrak{H}_2$, such that $A\in l^{\tau}_{\sigma_n}\mbox{form\,}A_1\triangle\, l^{\tau}_{\sigma_n}\mbox{form\,}A_2$.

On the other hand, by induction there are groups $B_1,\ldots ,B_a,$ $C_1,\ldots, C_b$, such that $B_k\in\mathfrak{F}_{i_k},$ $ C_k\in\mathfrak{F}_{j_k}$ and 
$$
A_1\in\nu_1(l^{\tau}_{\sigma_n}\mbox{form\,}B_1,\ldots,l^{\tau}_{\sigma_n}\mbox{form\,}B_n), \  A_2\in \nu_2(l^{\tau}_{\sigma_n}\mbox{form\,}C_1,\ldots,l^{\tau}_{\sigma_n}\mbox{form\,}C_b).
$$

Let the variables $x_{i_1},\ldots, x_{i_t}$ are not included in the word $\nu_2$, and all the variables $x_{i_{t+1}},\ldots, x_{i_a}$ are included in this word. Put $D_{i_k}=B_k$, if $k<t+1$, $D_{i_k}=B_k\times C_q$, where $q$ is such that $x_{i_k}=x_{j_q}$ for all $k\ge t+1$. Let $D_{j_k}=C_k$, if $x_{j_k}\notin\{x_{i_{t+1}},\ldots ,x_{i_b}\}$.
Let also $\mathfrak{L}_s=l^{\tau}_{\sigma_n}\mbox{form\,}D_{i_s}$, $\mathfrak{X}_r=l^{\tau}_{\sigma_n}\mbox{form\,}D_{j_r}$ for all $s=1,\ldots ,a$ and $r=1,\ldots ,b.$

Then $A_1\in\nu_1(\mathfrak{L}_1,\ldots ,\mathfrak{L}_a),$ $A_2\in\nu_2(\mathfrak{X}_1,\ldots ,\mathfrak{X}_b).$ Hence, there are formations $\mathfrak{H}_1,\ldots ,\mathfrak{H}_m$, such that
$$
A\in\nu_1(\mathfrak{H}_{i_1},\ldots ,\mathfrak{H}_{i_a})\triangle\,\nu_2(\mathfrak{H}_{j_1},\ldots ,\mathfrak{H}_{j_b})=\nu (\mathfrak{H}_1,\ldots ,\mathfrak{H}_m),
$$
where $\mathfrak{H}_j=l^{\tau}_{\sigma_n}\mbox{form\,}K_j$ and $K_j\in\mathfrak{F}_j$, $j=1,\ldots , m$.

Thus, the lattice $l^{\tau}_{\sigma_n}$ is $\mathfrak{G}$-separable. The theorem is proved.

In the case when $\tau$ is a trivial subgroup functor from Theorem~\ref{t11} we have 

\begin{corollary} {\rm \cite{voridr}}  \label{c37t11} 
The lattice $l^{\sigma}_n$ is $\mathfrak{G}$-separable.
\end{corollary}

In the classical case, when $\sigma=\sigma^1$ from Theorem~\ref{t11} we obtain the following well-known results

\begin{corollary} {\rm \cite[p. 159]{sk1}}  \label{c38t11} 
The lattice $l^\tau_n$ $\mathfrak{G}$ is separable.
\end{corollary}

\begin{corollary} {\rm \cite[p. 100]{shsk}} \label{c39t11} 
Let $A\in \omega(\mathfrak F_1,\ldots, \mathfrak F_m)$, where $\omega(x_1,\ldots , x_m)$ is the signature term $\{\wedge,\vee_n\}$, $\mathfrak{F}_1,\ldots ,\mathfrak{F}_m$ some $n$-multiply local formations. Then there are groups $A_1,\ldots, A_m$, $(A_i\in\mathfrak F_i)$, such that $A\in\omega(\mathfrak{M}_1,\ldots ,\mathfrak{M}_m)$, where $\mathfrak M_i=l_n\text{\rm form\,}A_i$.
\end{corollary}

\end{document}